\documentclass[11pt]{amsart}
\usepackage{amsmath,amssymb,epsfig,color}
\newtheorem{theorem}{Theorem}[section]

\newtheorem{corollary}[theorem]{Corollary}
\newtheorem{definition}[theorem]{Definition}
\newtheorem{lemma}[theorem]{Lemma}

\newtheorem{example}[theorem]{Example}

\newcommand{\vanish}[1]{}\parskip=10pt

\def\p{\prime}
\def\e{\varepsilon}
\def\B{\backslash}

\begin{document}

\title[Tensor Products of Signed Graphs]{Tutte Polynomials of Tensor Products of Signed Graphs
and their Applications in Knot Theory}
\author{Y. Diao, G. Hetyei and K. Hinson}
\address{Department of Mathematics and Statistics, UNC Charlotte,
    Charlotte, NC 28223}
\email{ydiao@uncc.edu, ghetyei@uncc.edu, kehinson@uncc.edu}
\dedicatory{} \subjclass{57M25}
\keywords{knots, Jones polynomials, Tutte polynomials, signed
graphs, tensor product of graphs.}
\begin{abstract}
It is well-known that the Jones polynomial of an alternating knot
is closely related to the Tutte polynomial of a special graph
obtained from a regular projection of the knot. Relying on the
results of Bollob\'as and Riordan, we introduce a generalization
of Kauffman's Tutte polynomial of signed graphs for which
describing the effect of taking a signed tensor product of signed
graphs is very simple. We show that this Tutte polynomial of a
signed tensor product of signed graphs may be expressed in terms
of the Tutte polynomials of the original signed graphs by using a
simple substitution rule. Our result enables us to compute the
Jones polynomials of some large non-alternating knots. The
combinatorics used to prove our main result is similar to Tutte's
original way of counting ``activities'' and specializes to a new,
perhaps simpler proof of the known formulas for the ordinary Tutte
polynomial of the tensor product of unsigned graphs or matroids.
\end{abstract}

\maketitle
\section{Introduction}

It is well-known that the Jones polynomial of a knot is related to
the Tutte polynomial of a special graph obtained from a regular
projection of the knot \cite{Ja,K1,K2}. It should be noted
however, that Tutte's original polynomial~\cite{Tu} is most easily
related to an {\em alternating} knot. If we want to compute the
Jones polynomial of a non-alternating knot, we need to use a
signed generalization of the Tutte polynomial, such as the one
introduced by Kauffman~\cite{K2}. Tutte defined his
polynomial~\cite{Tu} in terms of counting {\em activities} with
respect to a specific labelling of the graph, and his main result
is showing that the polynomial he introduced is independent of the
labelling, thus truly an invariant of the graph. This definition
and result may be easily generalized to matroids. The greatest
challenge in generalizing Tutte's polynomial to signed graphs or
matroids is to preserve the independence of the labelling. This
challenge is typically met by considering the Tutte polynomial of
a signed graph as an element of a polynomial ring modulo certain
relations between the variables. The most general result in the
area is due to Bollob\'as and Riordan~\cite{BR} who give a
necessary and sufficient set of relations modulo which a Tutte
polynomial of a signed graph is labelling independent. In
particular, Kauffman's Tutte polynomial for signed graphs is a
homomorphic image of the most general Tutte polynomial introduced
by Bollob\'as and Riordan.

Our main result consists of generalizing the known formulas for the
Tutte polynomial of the tensor product of two matroids to signed graphs
and matroids. The tensor product of matroids was introduced by
Brylawski~\cite{B} who also expressed the Tutte polynomial of a tensor
product in terms of the Tutte polynomials of the original
graphs~\cite{BO}. This result allowed Jaeger, Vertigan and
Welsh~\cite{Ja} to express the Jones polynomials of some alternating
knots. Our signed generalization allows us to compute the Jones polynomials
of some similarly complex knots that do not need to be alternating.

In the Preliminaries we review the Bollob\'as-Riordan definition of a
signed Tutte polynomial, using a notation that is closer to Kauffman's
paper~\cite{K2}. We introduce a signed Tutte polynomial in a factor ring
that is a proper homomorphic image of the most general signed Tutte
polynomial but, using a corollary from~\cite{BR} and a result from the
theory of determinantal rings, we observe that the signed Tutte
polynomial we propose is a most general possible, if we consider only
factor rings that are integral domains. In particular, our proposed
signed Tutte polynomial still allows us to compute Kauffman's signed Tutte
polynomial and later the Jones polynomial of an associated knot.

We introduce our definition of a signed tensor product of signed graphs
in Section~\ref{s_tp}. Since our aim is to apply this notion to graphs
associated to knots, our terminology implies considering connected
graphs only. However, our definitions, statements and proofs may be
generalized without any substantial change to signed matroids, very much
by the same reasons as the ones observed in~\cite[Remark 3]{BR}.

In Section~\ref{s_tcl} we generalize two polynomials associated to
a graph with a distinguished edge, used by Jaeger, Vertigan and
Welsh~\cite{Ja} in their formula for the unsigned tensor product.
These polynomials are defined by a system of equations
in~\cite{Ja}, we define their generalizations in terms of
(labelling dependent) ``activities'' in the spirit of Tutte's
paper~\cite{Tu}, and then show that our polynomials satisfy a
system of equations generalizing the defining equations
of~\cite{Ja}. The fact that these polynomials we introduce are
labelling independent follows from the equations having a unique
solution in an integral domain. This section is the reason why we
work with a most general integral domain but not the most general
ring in which a Tutte polynomial of signed graphs may be defined.

Our main result is in Section~\ref{s_tst}. We provide an explicit
substitution rule expressing the Tutte polynomial of a signed
tensor product of signed graphs in terms of the Tutte polynomials
of the original graphs. After showing that the proposed
substitution induces an endomorphism of the ring to which our
Tutte polynomials belong, it is sufficient to verify our rule for
one specific set of representatives, associated to one specific
labelling of the edges. The proof then becomes a classic
enumeration of ``activities'' the way, we believe, Tutte would
have proved this result. Our proof may also be applied to unsigned
graphs, and is perhaps more accessible to a non-expert than the
known proofs of the tensor product formulas for unsigned graphs
and matroids.

Our study is partly motivated by applications of knot theory to
physical and biological sciences in recent years, where there is a
need to distinguish different knots and measure the complexity of
knots. One common measure of knot complexity is the crossing
number of a knot, which is defined as the minimum number of
crossings over all possible regular projections of all possible
space realizations of the given knot. A remarkable result derived
from the Jones polynomial states that if a knot has a reduced
alternating projection (i.e. the knot is an alternating knot),
then the crossing number of the knot is equal to the number of
crossings in this reduced alternating projection \cite{K1,Mu,T0}.
Furthermore, for non-alternating knots, one can also use the Jones
polynomial to approximate their crossing numbers since it is known
that the span (or breadth) of the Jones polynomial of a knot
bounds the crossing number of the knot from below. Unfortunately,
the computation of the Jones polynomial of a knot is known to be
NP-hard \cite{Ja}. This prevents the computation of the Jones
polynomial of knots with large crossing numbers in general. For
example, \textit{Knotscape}, a commonly used software (developed
by Hoste and Thistlethwaite \cite{knotscape}) for computing
various knot invariants, handles knot diagrams up to about 50
crossings only. For some (non-alternating) knots constructed using
the tensor products of signed graphs, it is possible to compute
the Tutte polynomials of their corresponding graph tensor products
easily, leading to a quicker computation of their Jones
polynomials, hence a better understanding of such knots. Some such
examples are provided in Section~\ref{s_app}.

\section{Preliminaries}
\label{s_p}

In this section, we introduce the concept of a Tutte polynomial
for a signed graph, following closely the definitions in
\cite{BR}. The Tutte polynomial defined in this paper is only
slightly less general than the one defined in \cite{BR} and the
generalized Kauffman bracket polynomial \cite{K1,K2} may still be
obtained from it by factoring.

\begin{definition}
Let $G$ be a graph with edges labelled $1, 2, \ldots, n$, and let
$T$ be a spanning tree of $G$. An edge $e$ of $T$ is said to be
internally active if for any edge $f\not=e$ in G such that $(T
\setminus  e ) \cup  f $ is a spanning tree of $G$, the label of
$e$ is less than the label of $f$. Otherwise $e$ is said to be
internally inactive. On the other hand, an edge $f$ of $G
\setminus T$ is said to be externally active if $f$ has the
smallest label among the edges in the unique cycle contained in $T
\cup   f  $. Otherwise, $f$ is said to be externally inactive.
\end{definition}

\begin{example}\label{example1}{\em
Figure \ref{fig1} shows a graph with six edges and four vertices
where a labelling of the edges is given. The edges of a spanning
tree $T$ are highlighted in the graph. With respect to the tree
$T$, edges 1, 4 and 5 are internal edges and edges 2, 3 and 6 are
external edges. Edge 1 is internally active since it has the
smallest label, edge 4 is internally inactive since 4 is larger
than 3 in the spanning tree 315. Edge 5 is also internally
inactive since 5 is larger than 2 in the spanning tree 142. Edge 2
is externally active since 2 is the smallest in the cycle 245.
Edge 3 is externally inactive since 3 is larger than 1 in the
cycle 134 and edge 6 is apparently externally inactive since 6 is
the largest of all the labels.
\begin{figure}[!htb]
\begin{center}
\includegraphics[scale=0.8]{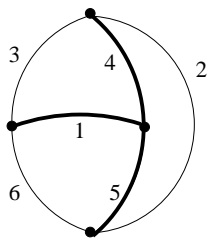}
\caption{An example of labelled graph with a marked spanning
tree.}\label{fig1}
\end{center}
\end{figure}}
\end{example}

\vspace{-0.4cm} Let $G$ be a signed (and connected) graph, that
is, each edge in $G$ is assigned a $+$ or $-$ sign. Let $T$ be a
spanning tree of $G$, for each edge $e$ in $G$ we will then assign
one of the following variables to it according to the activities
of $e$ (with respect to the tree $T$):

\begin{table}[h]\label{table1}
\begin{center}
\begin{tabular}{|c|c|c|}\hline
sign of $e$ & activity & variable assignment\\\hline\hline $+$ &
internally active & $x_+$\\\hline $-$ & internally active &
$x_-$\\\hline $+$ & externally active & $y_+$\\\hline $-$ &
externally active & $y_-$\\\hline $+$ & internally inactive &
$A_+$\\\hline $-$ & internally inactive & $A_-$\\\hline $+$ &
externally inactive & $B_+$\\\hline $-$ & externally inactive &
$B_-$\\\hline
\end{tabular}
\end{center}

\vspace{0.3cm} \caption{The variable assignment of an edge with
respect to a spanning tree $T$.}
\end{table}

\begin{definition} Let $G$ be a connected signed graph. For a
spanning tree $T$ of $G$, let $C(T)$ be the product of the
variable contributions from each edge of $G$ according to the
variable assignment above, then the Tutte polynomial $T(G)$ is
defined as the sum of all the $C(T)$'s over all possible spanning
trees of $G$.
\end{definition}

\begin{example}\label{example2}{\em
The graph given in Example \ref{example1} is assigned the signs as
shown in Figure \ref{fig2} below.
\begin{figure}[!htb]
\begin{center}
\includegraphics[scale=0.8]{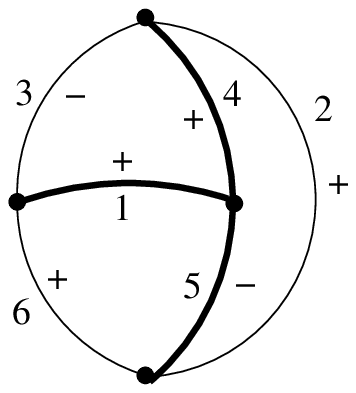}
\caption{An example of signed and labelled graph with a marked
spanning tree.}\label{fig2}
\end{center}
\end{figure}

For the spanning tree shown in the figure, the total contributions
of all the edges is easily calculated to be $x_+y_+A_+A_-B_+B_-$.
The complete list of the spanning trees of $G$ is given in Figure
\ref{fig3}. We leave it for our reader to verify that the Tutte
polynomial of $G$ is
\begin{eqnarray*}
& &4x_+y_+A_+A_-B_+B_- + 2y_+^{2}A_+A_-^{2}B_+ +
2x_+^{2}A_+B_+B_-^{2}\\
& +& y_+^{2}A_+^{2}A_-B_- + x_+^{2}A_-B_+^{2}B_- +
y_+^{2}y_-A_+^{2}A_- + x_+^{2}x_-B_+^{2}B_- \\
&+& y_+A_+A_-^{2}B_+^{2} + y_+A_+^{3}B_-^{2} + x_+A_-^{2}B_+^{3} +
x_+A_+^{2}B_+B_-^{2}.
\end{eqnarray*}

\begin{figure}[!htb]
\begin{center}
\includegraphics[scale=0.6]{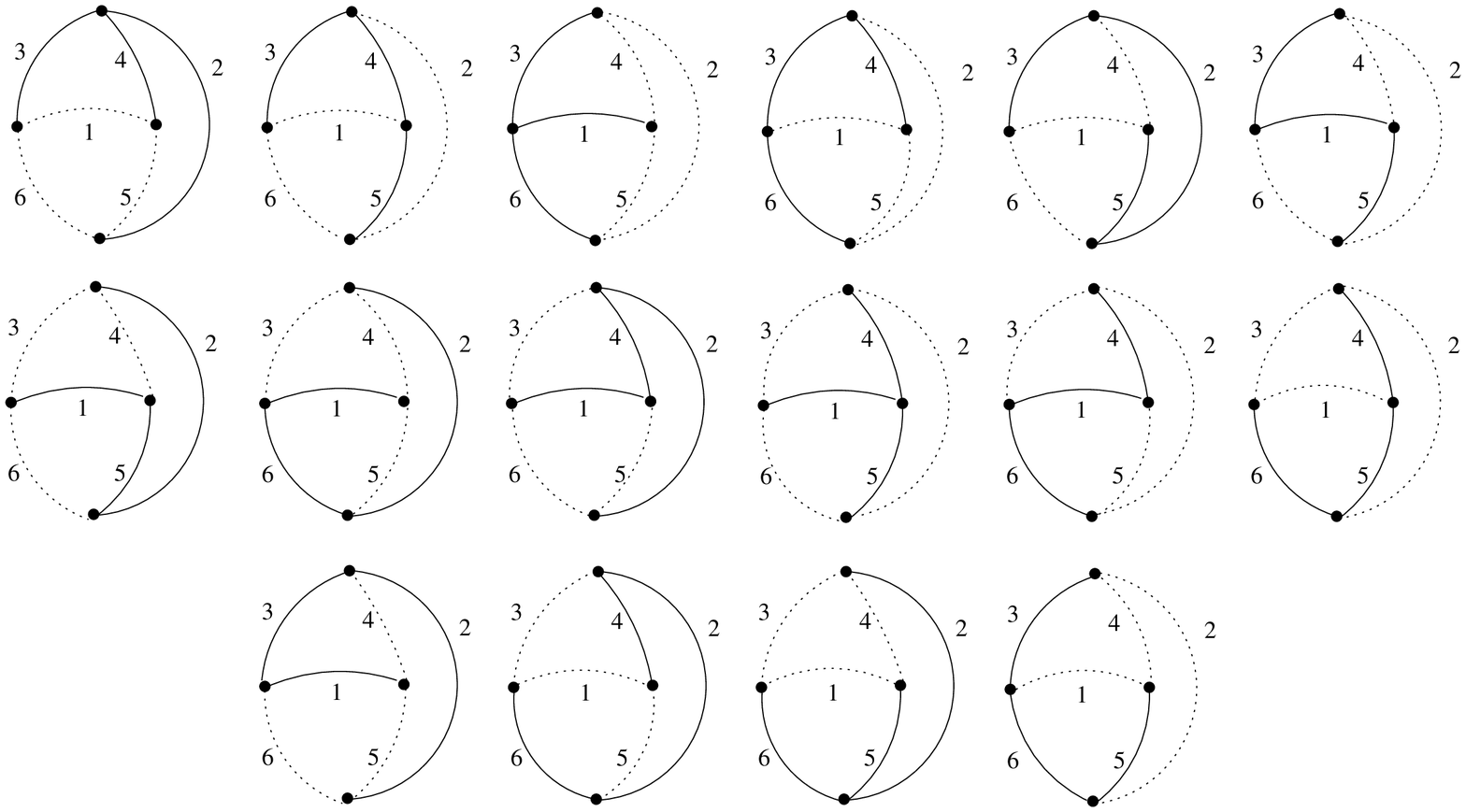}
\caption{The complete set of spanning trees of the graph given in
Figure \ref{fig1}.}\label{fig3}
\end{center}
\end{figure}
}
\end{example}

Notice that the Tutte polynomial defined this way is labelling
dependent. To remedy the situation, we want to factor the
polynomial ring ${\mathbb Z}[\Lambda]:={\mathbb
Z}[A_+,A_-,B_+,B_-,x_+,x_-,y_+,y_-]$ with an appropriate ideal
$I$, such that the formula for $T(G)$ in ${\mathbb Z}[\Lambda]/I$
becomes labelling independent. An exact description of all such
ideals (for a larger class of signed graphs) was given by
Bollob\'as and Riordan~\cite[Theorem 2]{BR}. (Here we state the
two-colored version.)
\begin{theorem}[Bollob\'as-Riordan]
\label{T_BR} Assume $I$ is an ideal of ${\mathbb Z}[\Lambda]$.
Then the homomorphic image of $T(G)$ in ${\mathbb Z}[\Lambda]/I$
is independent of the labelling of the edges of $G$ if and only if
\begin{equation}
\det\left(\begin{array}{ll}
x_+& B_+\\
x_-& B_-
\end{array}\right)
- \det\left(\begin{array}{ll}
A_+& y_+\\
A_-& y_-
\end{array}\right)\in I,
\end{equation}
\begin{equation}
y_{\nu}\det\left(\begin{array}{ll}
A_+& y_+\\
A_-& y_-
\end{array}\right)
- y_{\nu}\det\left(\begin{array}{ll}
A_+& B_+\\
A_-& B_-
\end{array}\right)\in I \quad\mbox{for $\nu\in\{+,-\}$,}
\end{equation}
and
\begin{equation}
x_{\nu}\det\left(\begin{array}{ll}
A_+& y_+\\
A_-& y_-
\end{array}\right)
- x_{\nu}\det\left(\begin{array}{ll}
A_+& B_+\\
A_-& B_-
\end{array}\right)\in I\quad\mbox{for $\nu\in\{+,-\}$.}
\end{equation}
\end{theorem}

Bollob\'as and Riordan denote the ideal generated by the
differences listed in Theorem~\ref{T_BR} by $I_0$. The homomorphic
image of $T(G)$ in ${\mathbb Z}[\Lambda]/I_0$ is the most general
signed Tutte polynomial whose definition is independent of the
labelling. We will describe how to factor ${\mathbb Z}[\Lambda]$
by an ideal properly containing $I_0$ to get the Jones polynomial.
To simplify our calculations, we want to replace $I_0$ with a
larger ideal $I_1$ in such a way that the effect of certain
operations of signed graphs is still describable in terms of Tutte
polynomials (as elements of ${\mathbb Z}[\Lambda]/I_1$). We will
keep our calculations as simple as possible since we want to
obtain a homomorphic image in an integral domain, and send $A_+$,
$A_-$, $B_+$, $B_-$, $x_+$, $x_-$, $y_+$, and $y_-$ into nonzero
entries at the end. Inspired by \cite[Corollary 3]{BR} we make the
following definition.
\begin{definition}\label{def2.6}
\label{D_oT} We consider the signed Tutte polynomial to be an
element of ${\mathbb Z}[\Lambda]/I_1$ where $I_1$ is the ideal
generated by
\begin{equation}
\label{E_I11}
\det\left(\begin{array}{ll}
x_+& B_+\\
x_-& B_-
\end{array}\right)
- \det\left(\begin{array}{ll}
A_+& B_+\\
A_-& B_-
\end{array}\right)
\end{equation}
and
\begin{equation}
\label{E_I12}
\det\left(\begin{array}{ll}
A_+& B_+\\
A_-& B_-
\end{array}\right)
- \det\left(\begin{array}{ll}
A_+& y_+\\
A_-& y_-
\end{array}\right).
\end{equation}
\end{definition}

Clearly $I_1$ properly contains $I_0$ so our Tutte polynomial is
labelling independent. Moreover, we highlight the following observation,
making \cite[Corollary 3]{BR} truly useful.

\begin{lemma}
\label{L_ip}
The ideal $I_1$ is a prime ideal. More generally, given any integral
domain $R$, the ideal $I_1$ generated by the elements
(\ref{E_I11}) and (\ref{E_I12}) in $R[\Lambda]$ is prime.
\end{lemma}
\begin{proof}
The linear map induced by $x_+\mapsto \widetilde{x}_{+}:=x_+-A_+$,
$x_-\mapsto\widetilde{x}_{-}:=x_--A_-$,
$y_+\mapsto\widetilde{y}_{+}:=y_+-B_+$,
$y_-\mapsto\widetilde{y}_{-}:=y_+-A_+$, $A_+\mapsto A_+$, $A_-\mapsto
A_-$, $B_+\mapsto B_+$ and $B_-\mapsto B_-$ is an isomorphism between
the polynomial rings $R[\Lambda]$ and
$R[\widetilde{\Lambda}]:=R[A_+,A_-,B_+,B_-,
  \widetilde{x}_{+},\widetilde{x}_{-},\widetilde{y}_{+},
  \widetilde{y}_{-}]$. Under this isomorphism, $I_1$ goes into the ideal
$\widetilde{I_1}$  generated by
$$
\Delta_1:=\det\left(\begin{array}{ll}
\widetilde{x}_+& B_+\\
\widetilde{x}_-& B_-
\end{array}\right)
\quad\mbox{and}\quad
\Delta_2:=\det\left(\begin{array}{ll}
A_+& \widetilde{y}_+\\
A_-& \widetilde{y}_-
\end{array}\right).
$$
We only need to show that $\widetilde{I_1}$ is a prime ideal in
$R[\widetilde{\Lambda}]$ or, equivalently, that the factor ring
$R[\widetilde{\Lambda}]/\widetilde{I_1}$ is an integral domain. We
prove the second equivalent statement in two steps. First we
consider the polynomial ring
$R[\widetilde{x}_+,\widetilde{x}_-,B_+,B_-]$ and the ideal
$\Delta_1\cdot R[\widetilde{x}_+,\widetilde{x}_-,B_+,B_-]$ in it.
This ideal is a {\em determinantal ideal}, and the factor ring
$R':=R[\widetilde{x}_+,\widetilde{x}_-,B_+,B_-]/(\Delta_1\cdot
R[\widetilde{x}_+,\widetilde{x}_-,B_+,B_-])$  is a {\em
  determinantal ring} as defined in the book of Bruns and Vetter~\cite{BV}.
As a consequence of Theorem (2.10) in~\cite{BV}, given any
integral domain $R$ and a square matrix $X$ of variables, the
ideal generated by the determinant of $X$ in the polynomial ring
$R[X]$ is prime. Thus $\Delta_1\cdot
R[\widetilde{x}_+,\widetilde{x}_-,B_+,B_-]$ is a prime ideal and
$R'$ is an integral domain. Observe now that
$R[\widetilde{\Lambda}]/\widetilde{I_1}$ is isomorphic to the
factor of $R'[\widetilde{y}_+,\widetilde{y}_-,A_+,A_-]$ modulo the
ideal $\Delta_2\cdot R'[\widetilde{y}_+,\widetilde{y}_-,A_+,A_-]$.
Applying \cite[Theorem (2.10)]{BV} again, this time to
$R'[\widetilde{y}_+,\widetilde{y}_-,A_+,A_-]$ and $\Delta_2\cdot
R'[\widetilde{y}_+,\widetilde{y}_-,A_+,A_-]$, we obtain that
$\Delta_2\cdot\widetilde{B}[\widetilde{y}_+,\widetilde{y}_-,A_+,A_-]$
is a prime ideal and
$R'[\widetilde{y}_+,\widetilde{y}_-,A_+,A_-]/(\Delta_2\cdot
R'[\widetilde{y}_+,\widetilde{y}_-,A_+,A_-])$ is an integral
domain.
\end{proof}

Finally, since the Tutte polynomial is labelling independent under
our conditions, it is easy to see that we have the following
recursive formula
\begin{equation}\label{recur}
T(G)=B_\varepsilon T(G\backslash e)+A_\varepsilon T(G/e),
\end{equation}
where $e$ is any given edge of $G$ with sign $\varepsilon$,
$G\backslash e$ is the graph obtained from $G$ by deleting $e$ and
$G/e$ is the graph obtained from $G$ by contracting $e$.

Note that each spanning tree of $G$ either contains $e$ or does
not contain $e$.  If $e$ is in the tree, then the tree can be
identified with a spanning tree of $G/ e $.  If $e$ is not in the
tree, then the tree can be identified with a spanning tree of
$G\setminus  e $.

\section{The Tensor Product of Two Signed Graphs}
\label{s_tp}

We will define the tensor product of two signed graphs in this
section. Notice that tensor product is usually defined for
matroids (see \cite{BO}). Our definition of a tensor product also
defines the graph only up to the underlying matroid structure.
This should not represent a problem since the definition of Tutte
polynomials depends only on the underlying matroid structure. All
our definitions and statements may be generalized without
essential adjustment to matroids just as was noted in
\cite[Remark 3]{BR}. Replacing each occurrence of the term
``spanning tree'' with ``matroid basis'', and understanding the
words ``cycle'' and ``dependence'' in the matroid theoretic sense
provides the immediate generalization. In particular, we only use
the term ``spanning tree'' because we are interested in
knot-theoretic applications where the associated graphs are
connected most of the time, but our definitions and proofs work
the same way for spanning forests of disconnected graphs.

\begin{definition} Let $M$ and $N$ be two signed graphs.
The positive (negative) tensor product of $M$ and $N$, denoted by
$M\otimes_+N$ ($M\otimes_-N$), is the (signed) graph obtained by
replacing each positive (negative) edge of $M$ with a copy of
$N\backslash e$, where $e$ is a fixed edge of $N$ that is to be
identified with the replaced edge of $M$.
\end{definition}

\begin{figure}[!htb]
\begin{center}
\includegraphics[scale=0.8]{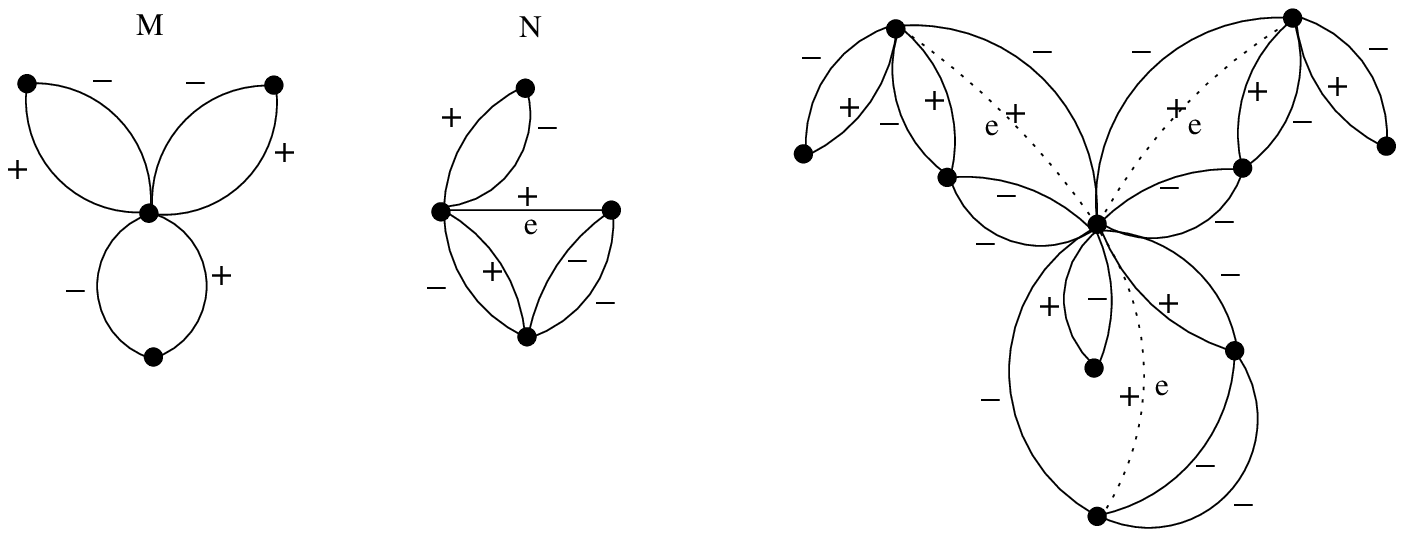}
\caption{The positive tensor product of two signed graphs $M$ and
$N$.}\label{figt}
\end{center}
\end{figure}

Figure \ref{figt} shows an example of a positive tensor product.
Notice that since there is no orientation defined on the edge $e$,
the two different ways of identifying $e$ with the replaced edge
in $M$ may lead to tensor products that are not isomorphic. Notice
also that the sign of the edge $e$ (in the graph $N$) does not
play any role in the tensor product since it does not appear in
the tensor product. In fact, an alternative definition of the
tensor product can be given by specifying two vertices of $N$
(that share the same boundary of a face of $N$) and replacing any
positive (negative) edge of $M$ with $N$ by identifying the two
vertices of the removed edge of $M$ with the two chosen vertices
on $N$.

For unsigned graphs (corresponding to alternating knots), formulas
for two special tensor products called the $k$-thickening and
$k$-stretch may be found in \cite[Lemma 6.3.24]{BO}, \cite[(7.2)
and (7.3)]{Ja}, and \cite[(3.8) and (3.10)]{Hu}. The general
formula for the general unsigned graphs can be found in \cite{BO}.
Our aim is to extend the known results in the unsigned cases to
the cases of signed graphs. Since our motivation is to apply these
results in knot theory, it makes sense to look at a few examples
that help one to understand what tensor products mean in knot
theory. The left side of Figure \ref{fignonalt} below shows a
signed graph $M$ with four vertices, at its right side is its
corresponding knot projection diagram. (The first paragraph of
Section \ref{s_app} discusses how to convert $M$ into its
corresponding knot projection diagram.)  Figure \ref{figrep1}
shows the graph $N$ with a distinguished edge $e$, as well as its
corresponding knot projection diagram. The tensor product
$M\otimes_+N$ is the same as $M$ (so the corresponding knot
projection does not change), while the tensor product
$M\otimes_-N$ is a graph with only positive edges whose
corresponding knot diagram is a reduced alternating knot
projection diagram as shown in Figure \ref{figalt}.

\begin{figure}[!htb]
\begin{center}
\includegraphics[scale=0.8]{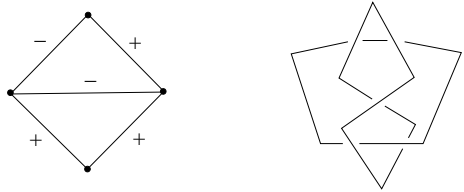}
\caption{A signed graph $M$ and its corresponding non-alternating
knot projection diagram.}\label{fignonalt}
\end{center}
\end{figure}

\begin{figure}[!htb]
\begin{center}
\includegraphics[scale=0.8]{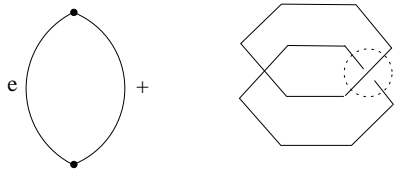}
\caption{A signed graph $N$ with a distinguished edge and its
corresponding knot projection diagram.}\label{figrep1}
\end{center}
\end{figure}

\begin{figure}[!htb]
\begin{center}
\includegraphics[scale=0.8]{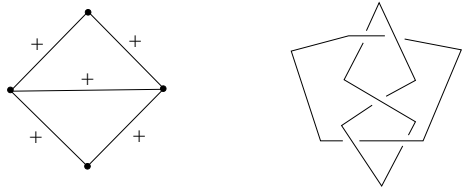}
\caption{The tensor product $M\otimes_-N$ and its corresponding
alternating knot projection diagram.}\label{figalt}
\end{center}
\end{figure}

This shows that at the simplest level, the tensor product
operation means changing the under/over strand information at the
crossings to make the new knot diagram an alternating one.

\section{The polynomials $T_C(N,e)$ and $T_L(N,e)$}
\label{s_tcl}

In this section we introduce two polynomials defined over a signed
graph $N$ with a distinguished edge $e$.  These polynomials are
signed generalizations of the polynomials $T_C$ and $T_L$ that
occur in the paper of Jaeger, Vertigan and Welsh~\cite{Ja}. A
variant of these unsigned polynomials was first introduced by
Brylawski~\cite{B}. The definition is usually done by a system of
equations. Here we make a definition that appears to be labelling
dependent, in terms of activities, in the spirit of Tutte's
original paper~\cite{Tu}. We will then show that our definition is
labelling independent. Our methods may be specialized to the
``traditional'' unsigned case, providing a new approach that
appears to be closer to Tutte's original way of thinking.

For any graph $N$ with a spanning tree $P$, any edge $f$ not on
$P$ will close a unique cycle in $P\cup f$. As usual in matroid
theory, we will use the notation $C(P,f)$ to denote this unique
cycle.  As mentioned after equation (6), $P$ can be identified
with a spanning tree of $N/ e $ or of $N\setminus  e $, depending
on whether $P$ contains $e$ or not. In the following, when dealing
with $N/ e $ we will use the shorthand of referring to a cycle
that contains the vertex contracted from $e$ as simply a cycle
containing $e$.

\begin{definition}
\label{D_tctl}
Let $N$ be a signed graph with a distinguished edge $e$.
Then $T_L(N,e)$ is the polynomial defined by the same rule that
defines $T(N\backslash e)$ except that internally active edges on
a cycle closed by $e$ will be considered as internally inactive
instead. Similarly, $T_C(N,e)$ is the polynomial defined by the
same rule that defines $T(N/e)$ except that externally active
edges that would close a cycle containing $e$ will be considered
as externally inactive instead.
\end{definition}

\begin{lemma}\label{LC}
Let $N$ be a graph, $P$ be a spanning tree of $N$, $e$ be an edge
on $P$, $e_i$ and $e_j$ be two distinct edges not on $P$ such that
$C(P,e_i)$ and $C(P,e_j)$ both contain the edge $e$. Then (1)
$P_i=(P\backslash e)\cup e_i$ and $P_j=(P\backslash e)\cup e_j$
are spanning trees of $N$ as well; (2) there exists at least one
cycle in $(P\backslash e)\cup e_i\cup e_j$; (3) any cycle in
$(P\backslash e)\cup e_i\cup e_j$ must contain both $e_i$ and
$e_j$.
\end{lemma}

\begin{proof}
(1) Since the only cycle in $P\cup e_i$ is $C(P,e_i)$ and $e$ is
on this cycle, removing $e$ from $P\cup e_i$ results in a graph
containing no cycles and having the same number of
edges as $P$, so $P_i$ is a spanning tree of $N$. Similarly, $P_j$ is also
a spanning tree of $N$.

(2) This is a classical graph theory result, which is equivalent to
one of the circuit axioms in matroid theory.

(3) Let $C$ be a cycle contained in $(P\backslash e)\cup e_i\cup
e_j$. If $C$ does not contain $e_i$, then it is totally contained
in $P_j$, but that is impossible since $P_j$ is a tree. So $C$
must contain $e_i$. Similarly, $C$ must contain $e_j$ as well.
\end{proof}

The following lemma is crucial in the proof of Theorem
\ref{gentctn}.

\begin{lemma}\label{bigequ}
The following identity holds in the ring ${\mathbb
Z}[\Lambda]/I_1$ for all $k\geq 1$ and all
$\varepsilon,\varepsilon_1,\ldots,\varepsilon_k\in\{-,+\}$.
$$
A_{\varepsilon}\left(\prod_{i=1}^k y_{\varepsilon_i} -\prod_{i=1}^k
B_{\varepsilon_i} \right)
=
(y_{\varepsilon}-B_{\varepsilon}) \sum_{i=1}^k A_{\varepsilon_i}
\prod_{j=1}^{i-1} y_{\varepsilon_j}
\prod_{t=i+1}^k B_{\varepsilon_t}.
$$
As usual, all empty products are equal to $1$.
\end{lemma}
\begin{proof}
We proceed by induction on $k$. Assume first $k=1$. For
$\varepsilon_1=\varepsilon$ the identity
$$
A_{\varepsilon}(y_{\varepsilon}-B_{\varepsilon})=
(y_{\varepsilon}-B_{\varepsilon})A_{\varepsilon}
$$
is valid even in ${\mathbb Z}[\Lambda]$, whereas for
$\varepsilon_1\neq\varepsilon$
the identity
$$
A_{\varepsilon}(y_{\varepsilon_1}-B_{\varepsilon_1})=
(y_{\varepsilon}-B_{\varepsilon})A_{\varepsilon_1}
$$
may be rearranged as
$$
\det\left(\begin{array}{ll}
A_{\varepsilon_1}& B_{\varepsilon_1}\\
A_{\varepsilon}& B_{\varepsilon}
\end{array}\right)
-
\det\left(\begin{array}{ll}
A_{\varepsilon_1}& y_{\varepsilon_1}\\
A_{\varepsilon}& y_{\varepsilon}
\end{array}\right)=0,
$$
and the left hand side is one of the generators of $I_1$.
To prove the induction step, let us restate the identity in the
following equivalent form.
\begin{equation}
\label{E_sid}
\det\left(\begin{array}{ll}
\sum\limits_{i=1}^k A_{\varepsilon_i} \prod\limits_{j=1}^{i-1}
y_{\varepsilon_j}
\prod\limits_{t=i+1}^k B_{\varepsilon_t}
& \prod\limits_{i=1}^k B_{\varepsilon_i} \\
A_{\varepsilon} & B_{\varepsilon}\\
\end{array}\right)
=
\det\left(\begin{array}{ll}
\sum\limits_{i=1}^k A_{\varepsilon_i} \prod\limits_{j=1}^{i-1}
y_{\varepsilon_j}
\prod\limits_{t=i+1}^k B_{\varepsilon_t} & \prod\limits_{i=1}^k
y_{\varepsilon_i}\\
A_{\varepsilon} & y_{\varepsilon}\\
\end{array}\right).
\end{equation}
Let us use $LHS_{k}$ and $RHS_{k}$ respectively as a short hand
for the left hand side and right hand side of (\ref{E_sid})
respectively. Assume by induction that (\ref{E_sid}) is valid for
some $k\geq 1$ and consider the identity obtained by increasing $k$ to
$k+1$. The new left hand side
$$
LHS_{k+1}=
\det\left(\begin{array}{ll}
\sum\limits_{i=1}^{k+1} A_{\varepsilon_i} \prod\limits_{j=1}^{i-1}
y_{\varepsilon_j}
\prod\limits_{t=i+1}^{k+1} B_{\varepsilon_t}
& \prod\limits_{i=1}^{k+1} B_{\varepsilon_i} \\
A_{\varepsilon} & B_{\varepsilon}\\
\end{array}\right)
$$
may be rewritten as
$$
LHS_{k+1}=
\det\left(\begin{array}{ll}
\left(\sum\limits_{i=1}^{k} A_{\varepsilon_i} \prod\limits_{j=1}^{i-1}
y_{\varepsilon_j}
\prod\limits_{t=i+1}^{k} B_{\varepsilon_t}\right)B_{\varepsilon_{k+1}}
+A_{\varepsilon_{k+1}} \prod\limits_{j=1}^{k} y_{\varepsilon_j}
& \prod\limits_{i=1}^{k} B_{\varepsilon_i}\cdot  B_{\varepsilon_{k+1}}\\
A_{\varepsilon} & B_{\varepsilon}\\
\end{array}\right),
$$
which yields the recursion
$$
LHS_{k+1}= LHS_{k}\cdot B_{\varepsilon_{k+1}}+A_{\varepsilon_{k+1}}
  \prod\limits_{j=1}^{k} y_{\varepsilon_j}\cdot B_{\varepsilon}.
$$
Similarly, the new right hand side may be rewritten as
$$
RHS_{k+1}=
\det\left(\begin{array}{ll}
\left(\sum\limits_{i=1}^{k} A_{\varepsilon_i} \prod\limits_{j=1}^{i-1}
y_{\varepsilon_j}
\prod\limits_{t=i+1}^{k} B_{\varepsilon_t}\right)B_{\varepsilon_{k+1}}
+A_{\varepsilon_{k+1}} \prod\limits_{j=1}^{k} y_{\varepsilon_j}
& \prod\limits_{i=1}^{k}y_{\varepsilon_i}\cdot y_{\varepsilon_{k+1}}\\
A_{\varepsilon} & y_{\varepsilon}\\
\end{array}\right),
$$
yielding the recursion
$$
RHS_{k+1}=RHS_k\cdot
B_{\varepsilon_{k+1}}
+\prod\limits_{i=1}^{k}y_{\varepsilon_i}
\left(A_{\varepsilon_{k+1}} y_{\varepsilon}
-A_{\varepsilon}\left(y_{\varepsilon_{k+1}}-B_{\varepsilon_{k+1}}\right)
\right).
$$
To show that the left hand side and the right hand side obey the same
recursion formula, it is sufficient to observe that
$$
A_{\varepsilon_{k+1}}B_{\varepsilon}=A_{\varepsilon_{k+1}} y_{\varepsilon}
-A_{\varepsilon}\left(y_{\varepsilon_{k+1}}-B_{\varepsilon_{k+1}}\right)
$$
holds in ${\mathbb Z}[\Lambda]/I_1$, which is equivalent to the already shown
induction basis statement.
\end{proof}

We will now state and prove our main result of this section.

\begin{theorem}\label{gentctn}
Let $N$ be a signed graph with a distinguished edge $e$ and
$\epsilon\in\{-1,+1\}$, then
\begin{equation}\label{TLC1}
A_\e (T(N/e)-T_C(N,e))=(y_\e-B_\e)T_L(N,e)
\end{equation}
and
\begin{equation}\label{TLC2}
B_\e(T(N\B e)-T_L(N,e))=(x_\e-A_\e)T_C(N,e).
\end{equation}
\end{theorem}

Note that the equalities in the above theorem can also be
expressed in terms of determinants:
\begin{equation}\label{tctnmat1}
\det\left(\begin{array}{ll}
T_L(N,e)& T_C(N,e)\\
A_\e & B_\e
\end{array}\right)
= \det\left(\begin{array}{ll}
T_L(N,e)& T(N/e)\\
A_\e& y_\e
\end{array}\right)
\end{equation}
and
\begin{equation}\label{tctnmat2}
\det\left(\begin{array}{ll}
T_L(N,e)& T_C(N,e)\\
A_\e& B_\e
\end{array}\right)
= \det\left(\begin{array}{ll}
T(N\backslash e)& T_C(N,e)\\
x_\e& B_\e
\end{array}\right)
\end{equation}

\begin{proof}
We will prove equation (\ref{TLC1}) first. Recall that $T_L(N,e)$
is defined by the same rule that defines $T(N\backslash e)$ except
that an internally active edge with sign $\e$ on a cycle closed by
$e$ is considered inactive, hence contributes $A_\e$ instead of
$x_\e$, and $T_C(N,e)$ is defined by the same rule that defines
$T(N/e)$ except that an externally active edge with sign $\e$ that
would close a cycle containing $e$ contributes $B_\e$ instead of
$y_\e$. By definition the only edges having different weight in
$T(N/e)$ and $T_C(N,e)$ are those externally active edges on a
cycle containing $e$ since the externally inactive edges are not
the exceptions in the definition of $T_C(N,e)$ and would make the
same contributions in $T(N/e)$ and $T_C(N,e)$. We will assign $e$
the largest index among the edges of $N$. Since $e$ does not
appear in $N/e$ nor in $N\B e$, this will not affect our
polynomials.

Let $Q$ be a spanning tree of $N$ with the property that $e\in Q$
and there exists at least one externally active edge $g$ with
respect to $Q$ such that $C(Q,g)$ contains $e$. Let $S$ be the set
of such $(Q,g)$ pairs. We will establish a one-to-one
correspondence between $S$ and the set $\mathcal{T}_1$ of spanning
trees of $N\B e$. Let $P$ be a spanning tree of $N\B e$. Since $e$
has the largest label, the edge $g$ with smallest label on the
cycle $C(P,e)$ is different from $e$. It follows that $Q=(P\B
g)\cup e$ is a spanning tree of $N$ and $g$ is an externally
active edge with respect to $Q$. This sets up the correspondence
between $P$ and the pair $(Q,g)$. We leave it to our reader to
verify that this is indeed a one-to-one correspondence.

On the other hand, the spanning trees of $N/e$ can be divided into
two disjoint sets: those that contain at least one externally active edge
that closes a cycle containing $e$ and those that do not. Denote
the first set by $\mathcal{T}_2$ and the second by $\mathcal{T}_3$. If
we pair a spanning tree in $\mathcal{T}_2$ and an externally active edge
that closes a cycle containing the vertex contracted from $e$,
then apparently this will set up a one-to-one correspondence with
the set $S$ as well. In other words, there is a one-to-one
correspondence between $\mathcal{T}_1$ and the set $\mathcal{T}_4$
of pairs of the spanning trees of $\mathcal{T}_2$ and their
corresponding edges described above.

Consider a spanning tree $P\in \mathcal{T}_3$. Since there are no
externally active edges with respect to $P$, the exception rule
defining $T_C(N,e)$ does not apply to any edge, hence the
contributions of $P$ to $T(N/e)$ and $T_C(N,e)$ are identical. It
follows that the combined contribution of $\mathcal{T}_3$ to the
LHS is $0$ and we may ignore this set entirely in the rest of our
argument.

Let $P^\p$ be a spanning tree of $N/e$ from the set
$\mathcal{T}_2$ and let $e_1$, $e_2$, ..., $e_k$ ($k\ge 1$) be the
externally active edges of $N/e$ with respect to $P^\p$ such that
$C(P^\p,e_j)$ contains the vertex contracted from $e$. Let $\e_j$
be the sign of $e_j$ ($1\le j\le k$) and let $P$ be the spanning
tree of $N$ obtained from $P^\p$ by recovering the edge $e$.
Notice that for each $j$, $(P\B e)\cup e_j$ is a spanning tree of
$N\B e$ by Lemma \ref{LC}. Furthermore,
$(P^\p,e_j)\longleftrightarrow (P\B e)\cup e_j=P_j$ is the
one-to-one correspondence between $\mathcal{T}_1$ and
$\mathcal{T}_4$ discussed above. The combined contribution of
$P^\p$ (or of all the pairs $(P^\p,e_j)$, $1\le j\le k$) to
$T(N/e)-T_C(N,e)$ is $(\prod_{j=1}^k y_{\e_j}-\prod_{j=1}^k
B_{\e_j})$ times the weight of the other edges of $P^\p$ according
to the exception rule of $T_C(N,e)$. Since (by Lemma \ref{bigequ})
$$
A_\e(\prod_{j=1}^k y_{\e_j}-\prod_{j=1}^k B_{\e_j})=(y_\e-B_\e)
\sum_{j=1}^kA_{\e_j}\prod_{i=1}^{j-1}y_{\e_i}\prod_{t=j+1}^{k}B_{\e_t},
$$
it suffices to show that, assuming that $e_1$, $e_2$, ..., $e_k$
are listed in increasing order of their labels, the contribution
of $P_j$ to $T_L(N,e)$ is
$A_{\e_j}\prod_{i=1}^{j-1}y_{\e_i}\prod_{t=j+1}^{k}B_{\e_t}$ times
the weight of the other edges of $P$ (not including $e$) for each
$j$, provided that the weight of the other edges on $P$ is the
same on the LHS and RHS of the equation in the theorem.

Case 1. The edges of the cycle $C_j=C(P,e_j)=C(P_j,e)$ that are
different from $e$ and $e_j$: These edges are internal in both
$P_j$ and $P$. On the side of $T(N/e)-T_C(N,e)$, these edges are
internally inactive due to the fact that $e_j$ is externally
active. On the side of $T_L(N,e)$, these edges are also treated as
inactive edges since they are exactly the exceptions in the
definition of $T_L(N,e)$. Thus, these edges of $P$ make the same
contributions on the two sides of the equation in the theorem.

Case 2. The edges of $P$ that are not contained in the cycle
$C_j$: Let $f$ be such an edge. Then $f$ is either active
in both $N/e$ (with respect to $P^\p$) and $N\B e$ (with
respect to $P_j$), or inactive in both,
because it is compared to the same edges in both
$N/e$ and $N\backslash e$, and none of the exceptions applies. In
fact, if $f\in C(P^\p,g)$ for some external edge $g$ then
$g\not=e_j$ (otherwise $f\in C_j$), so $g$ is an external edge
with respect to $P_j$ as well. The exception rule does not apply
to $f$ on the LHS since $f$ is internal and the exception rule
does not apply to $f$ on the RHS since $f\not\in C_j$. Conversely,
if $f\in C(P_j,g)$ for some $g$ external to $P_j$, then $g\not=e$
(otherwise $f\in C_j$ again). So $g$ is also an edge external with
respect to $P_j$. Thus, these edges of $P$ also make the same
contributions on the two sides of the equation in the theorem.

Case 3. The external edges of $N\backslash P$ with respect to
$P_j$ that are different from the $e_j$'s: Let $f$ be such an
external edge. If $f$ is externally inactive in $N\B e$ (with
respect to $P_j$), then it is also externally inactive in $N/e$
(with respect to $P^\p$). So its weight on both sides is the same.
If $f$ is externally active in $N\B e$ (with respect to $P_j$),
then $C(P,f)$ (hence $C(P^\p,f)$) does not contain the edge $e$
(or it would have to be one of the $e_j$'s). It follows that
$C(P^\p,f)=C(P_j,f)$. So $f$ is externally active in $N/e$ with
respect to $P^\p$ as well. Thus $f$ makes the same contribution on
both sides of the equation.

Case 4. We will now consider the weight of the edges $e_1$, $e_2$,
..., $e_{j+1}$, ..., $e_k$ in $T_L(N,e)$. Consider the unique
cycle $C(P_j,e_i)$ for any $i\not=j$. By Lemma \ref{LC}, $e_i$ and
$e_j$ are both contained in this cycle. It follows that $e_i$ is
externally inactive if $i>j$ by our assumption on the labels of
the edges $e_1$, $e_2$, ..., $e_k$. On the other hand, notice that
$C(P_j,e_i)$ is a subset of $C(P,e_i)\cup C(P,e_j)$. By the
definition of $e_i$ and $e_j$, $e_i$ has the smallest label among
the edges of $C(P,e_i)$ and $e_j$ has the smallest label among the
edges of $C(P,e_j)$. Thus, if $i<j$, then $e_i$ has the smallest
label among all edges of $C(P,e_i)\cup C(P,e_j)$. In particular,
$e_i$ has the smallest label among all edges of $C(P_j,e_i)$. That
is, $e_i$ is externally active with respect to $P_j$ (in $N\B e$).
So the contribution of $e_i$ is $y_{\e_i}$ for any $1\le i<j$. It
follows that the total contribution of the edges $e_1$, $e_2$,
..., $e_{j+1}$, ..., $e_k$ in $T_L(N,e)$ is
$\prod_{i=1}^{j-1}y_{\e_i}\prod_{t=j+1}^{k}B_{\e_t}$.

Case 5. Since $e_j$ is an internal edge with respect to $P_j$ and
it is on a cycle closed by $e$ (by the definition of $P_j$), $e_j$
is considered to be inactive by the exception rule of $T_L(N,e)$.
Thus the contribution of $e_j$ is $A_{\e_j}$.

Combining the above cases finishes our proof for (\ref{TLC1}). The
proof for (\ref{TLC2}) is quite similar to that of (\ref{TLC1}).
However, for the sake of convenience of our reader, we still
provide the details below.

$T_L(N,e)$ is defined by the same rule as $T(N\backslash
 e )$ except that internally active edges on a cycle closed by
$e$ are considered as inactive instead.  Thus the only edges
having different weight in $T(N\backslash  e )$ and $T_L(N,e)$ are
those internally active edges on a cycle closed by $e$, since the
internally inactive edges make the same contributions in
$T(N\backslash  e )$ and $T_L(N,e)$. Keep in mind that we have
assigned $e$ the largest index among the edges of $N$.

Let $P$ be a spanning tree in $N\backslash  e $ and let $e_1$,
$e_2$, ..., $e_k$ ($k\ge 1$) be the internally active edges such
that $P \backslash  e_j  \cup  e $ is a spanning tree of $N$. Let
$P^\p_j$ be the spanning tree of $N/ e $ obtained from $P$ by
contracting the edge $e$ and deleting the edge $e_j$ from the
cycle that is thus formed. Notice that $P^\p_j \longleftrightarrow
P\backslash e_j \cup e$ is a one-to-one correspondence. Since the
edges $e_1$, $e_2$, ..., $e_k$ are on the cycle closed by $e$ in
$P$, they are subject to the exception rule in the definition of
$T_L(N,e)$, and all other edges carry the same weight in the
calculations of $T_L(N,e)$ and $T(N\backslash e)$. Thus the
contribution of $P$ to $T(N\B e)-T_L(N,e)$ is
$(\prod_{j=1}^kx_{\e_j}-\prod_{j=1}^kA_{\e_j})$ times the combined
weight of the other edges.  Since
$$
B_\e(\prod_{j=1}^kx_{\e_j}-\prod_{j=1}^kA_{\e_j})=(x_\e-A_\e)
\sum_{j=1}^kB_{\e_j}\prod_{i=1}^{j-1}x_{\e_i}\prod_{t=j+1}^{k}A_{\e_t}
$$
by Lemma \ref{bigequ}, it suffices to show that, assuming that
$e_1$, $e_2$, ..., $e_k$ are listed in increasing order of their
labels, the contribution of $P^\p_j$ to $T_C(N,e)$ is
$B_{\e_j}\prod_{i=1}^{j-1}x_{\e_i}\prod_{t=j+1}^{k}A_{\e_t}$ times
the weight of the other edges of $P^\p_j$ for each $j$, provided
that the weight of the other edges on $P$ is the same on the LHS
and the RHS of equation (\ref{TLC2}) in the theorem.

If $g$ is an edge of $P$, let us denote by $K(P,g)$ the set of
external edges $f$ such that $(P \setminus  g ) \cup  f $ is a
spanning tree. We have the following cases to consider.

Case 1.  The edges of $K(P,e_j)$ that are different from $e$ and
$e_j$. These edges are external in both $P$ and $P\backslash e_j
\cup e$.  On the side of $T(N\backslash e) - T_L(N,e)$, these
edges are externally inactive as $e_j$ is internally active. On
the side of $T_C(N,e)$, these edges are the exceptions in the
definition of $T_C(N,e)$ hence are also considered as inactive.
Hence these edges make the same contributions on both sides of the
equation.

Case 2.  The edges of $N\backslash P$ that are not in $K(P,e_j)$.
Let $f$ be such an edge. Then $f$ is either inactive in both
$N\backslash e$ (with respect to $P$) and $N/e$ (with respect to
$P^\p_j$) or active in both $N\backslash e$ and $N/e$, because it
is compared to the same edges in both $N\backslash e$ and $N/e$.
The exception rule does not apply since $f$ is external on the LHS
and $e \not\in C(P^\p_j,f)$ on the RHS . Hence $f$ makes the same
contribution on both sides of the equation.

Case 3.  The internal edges of $P$ that are different from the
$e_i$'s. Let $f$ be such an edge. If $f$ is internally inactive in
$N/e$ (with respect to $P^\p_j$) then it is also internally
inactive in $N \backslash e$ (with respect to $P$). If $f$ is
internally active in $N/e$ then $K(P^\p_j, f)$ does not contain
the vertex contracted from $e$ (otherwise it would be one of the
$e_i$'s).  It follows that $K(P,f) = K(P\backslash  e_j  \cup e,
f)$, so $f$ is internally active with respect to $P$ as well. Thus
$f$ makes the same contribution on both sides of the equation.

Case 4.  We now consider the edges $e_1$, $e_2$, ..., $e_k$ in
$T_C(N,e)$. Recall that $P^\p_j$ is the spanning tree of $N/e$
obtained from $P$ by contracting the edge $e$ and deleting the
edge $e_j$. Consider the cut $K(P^\p_j, e_i)$ for any $i\not=j$.
$e_i$ and $e_j$ are both contained in this cut. It follows that
$e_i$ is internally inactive if $i>j$, by our assumption on the
labels of the edges $e_1$, $e_2$, ..., $e_k$, and that $e_i$ is
internally active if $i<j$.  So the contribution of $e_i$ is
$x_{\e_i}$ for any $1\le i<j$ and $A_{\e_i}$ for $j<i\le k$.
Finally, since $e_j$ is external to $P^\p_j$ and it closes a cycle
in $P^\p_j$ containing $e$, the exception rule for $T_C(N,e)$
applies to it. So it is always considered as inactive hence its
contribution in $T_C(N,e)$ is $B_{\e_j}$. It follows that the
total contribution of $e_1$, $e_2$, ..., $e_k$ in $T_C(N,e)$ is
$B_{\e_j}\prod_{i=1}^{j-1}x_{\e_i}\prod_{t=j+1}^k A_{\e_t}$, as
desired.
\end{proof}
One of the most important consequences of Theorem~\ref{gentctn} is the
following.
\begin{corollary}
The polynomials $T_C(N,e)$ and $T_L(N,e)$ are {\em independent of the
labelling}. They may be equivalently defined by the system of equations
(\ref{TLC1}) and (\ref{TLC2}).
\end{corollary}
In fact, by Theorem~\ref{gentctn}, for any $\e\in\{+,-\}$, setting
  $z_C=T_C(N,e)$ and $z_L=T_L(N,e)$ provides a solution of the linear
  system of equations
$$
\begin{array}{rcl}
(y_\e-B_\e)z_L+A_\e z_C&=&A_\e T(N/e)\\
B_\e z_L+(x_\e-A_\e) z_C&=&B_\e T(N\B e).\\
\end{array}
$$
Here $z_L$ are $z_C$ are the unknowns, and the givens belong to
${\mathbb Z}[\Lambda]/I_1$, an integral domain by
Lemma~\ref{L_ip}. Cramer's rule is applicable in the quotient
field of ${\mathbb
  Z}[\Lambda]/I_1$, and we have
$$
\det
\left(
\begin{array}{ll}
y_\e-B_\e & A_\e\\
B_\e & x_\e-A_\e\\
\end{array}
\right)=x_\e y_\e-A_\e y_\e-B_\e x_\e
$$
which is a nonzero element of ${\mathbb Z}[\Lambda]/I_1$, independently
of the choice of $\e$, since each element of $I_1$ is a ${\mathbb
  Z}$-linear combination of monomials involving both negative and
positive variables.

We may recover the classical unsigned case that is widely studied in the
literature in two steps as follows. Consider unsigned graphs
as special signed graphs whose edges are all positive. The signed Tutte
polynomial of such a graph may be considered as a polynomial in the
positive variables only. The subring of ${\mathbb Z}[\Lambda]/I_1$,
generated by the integers and the positive variables is easily seen to
be isomorphic to the polynomial ring ${\mathbb Z}[x_+,y_+,A_+,B_+]$,
because $I_1$ involves only polynomials with terms of mixed signature, as
mentioned above. We thus consider $T_C(N,e)$ and $T_L(N,e)$ as elements of
${\mathbb Z}[x_+,y_+,A_+,B_+]$, defined by the exception rules given
in Definition~\ref{D_tctl}. Under these circumstances Lemma~\ref{bigequ} may be
replaced with the trivial identity
$$
A_+(y_+^k-B_+^k)=A_+(y_+-B_+)\cdot \sum_{i=1}^k y_+^{i-1}B_+^{k-i},
$$
and Theorem~\ref{gentctn} specializes to the following statement.
\begin{theorem}
\label{utctn}
Given an unsigned graph $N$ with a distinguished edge $e$, the signed
polynomials  $T_C(N,e), T_L(N,e)\in {\mathbb Z}[x_+,y_+,A_+,B_+]$
satisfy
$$A_+ (T(N/e)-T_C(N,e))=(y_+-B_+)T_L(N,e)$$
and
$$
B_+(T(N\B e)-T_L(N,e))=(x_+-A_+)T_C(N,e).
$$
\end{theorem}
The system of equations in Theorem~\ref{utctn} uniquely determines
the polynomials $T_C(N,e), T_L(N,e)\in {\mathbb Z}[x_+,y_+,A_+,B_+]$,
and may be used as their alternative definition. Consider finally
the ring homomorphism $\phi: {\mathbb Z}[x_+,y_+,A_+,B_+]\rightarrow {\mathbb
  Z}[x,y]$, given by $x_+\mapsto x$, $y_+\mapsto y$, $A_+\mapsto 1$, and
$B_+\mapsto 1$. On the one hand  $T_C(N,e)$, $T_L(N,e)$ go into
polynomials that may be defined by modifying the definition of the
ordinary Tutte polynomials $T(N\B e;x,y)$ and $T(N/e;x,y)$
according to the exception rules given in Definition~\ref{D_tctl}.
On the other hand these homomorphic images satisfy the homomorphic
images of the equations given in Theorem~\ref{utctn} which still
has a unique solution in the polynomial ring ${\mathbb Z}[x,y]$.

\begin{corollary}
Given an unsigned graph $N$ with a distinguished edge
$e$, let us define the polynomials $T_C(N,e)$ and $T_L(N,e)$ by
modifying the definition of the ordinary Tutte polynomials $T(N\B
e;x,y)$ and $T(N/e;x,y)$
according to the exception rules given in Definition~\ref{D_tctl}. Then these
polynomials may be equivalently defined by the system of equations
\begin{equation}
\label{E_utctl}
\begin{array}{rcl}
T(N/e)-T_C(N,e)&=&(y-1)T_L(N,e)\\
T(N\B e)-T_L(N,e)&=&(x-1)T_C(N,e).\\
\end{array}
\end{equation}
In particular, the definition is labelling independent.
\end{corollary}
The system of equations (\ref{E_utctl}) appears as the definition of the
unsigned polynomials $T_C$ and $T_L$ in~\cite[Equation (4.2)]{Ja}.

\section{The Tutte polynomial of a signed tensor product}
\label{s_tst}

The main result of this section is that the Tutte polynomial of
the tensor product $M\otimes_+N$ or $M\otimes_-N$ can be computed
using the Tutte polynomial $T(M)$ of $M$ and the special Tutte
polynomials $T(N\B e)$, $T(N/e)$, $T_C(N,e)$ and $T_L(N,e)$ of $N$
discussed in the last section by simple variable substitutions.
First we need the following lemma.

\begin{lemma}\label{Lsub}
Let $N$ be a signed graph with a distinguished edge $e$. Then the
endomorphism of ${\mathbb Z}[\Lambda]$ given by
$$
x_+\mapsto T(N\backslash e)\quad A_+\mapsto T_L(N,e)\quad
y_+\mapsto T(N/e)\quad B_+\mapsto T_C(N,e)
$$
(with all negative variables unchanged) sends $I_1$ into itself.
Similarly, the endomorphism of ${\mathbb Z}[\Lambda]$ given by
$$
x_-\mapsto T(N\backslash e)\quad A_-\mapsto T_L(N,e)\quad
y_-\mapsto T(N/e)\quad B_-\mapsto T_C(N,e)
$$
(with all positive variables unchanged) sends $I_1$ into itself as
well. Consequently, these endomorphisms induce endomorphisms of
the factor ring ${\mathbb Z}[\Lambda]/I_1$.
\end{lemma}

\begin{proof}
This is immediate from Equations
(\ref{tctnmat1}) and (\ref{tctnmat2}), showing that the given endomorphisms
send the two generators of $I_1$ into elements of $I_1$.
\end{proof}

\begin{theorem}\label{Ttensor}
Let $M$ be a signed graph and $N$ a signed graph with a
distinguished edge $e$. Then $T(M\otimes_+N)$ can be computed from
$T(M)$ by keeping the negative variables unchanged and using the
substitutions
$$
x_+\mapsto T(N\backslash e)\quad
A_+\mapsto T_L(N,e)\quad
y_+\mapsto T(N/e)\quad
B_+\mapsto T_C(N,e).
$$
Similarly, $T(M\otimes_-N)$ can be computed from $T(M)$ by keeping
the positive variables unchanged and using the substitutions
$$
x_-\mapsto T(N\backslash e)\quad
A_-\mapsto T_L(N,e)\quad
y_-\mapsto T(N/e)\quad
B_-\mapsto T_C(N,e).
$$
\end{theorem}

\begin{proof}
By Lemma \ref{Lsub}, the given substitutions induce an
endomorphism of ${\mathbb Z}[\Lambda]/I_1$. This means that the
operation is well-defined in the sense that no matter which
representative of $T(M)$, $T(N\backslash e)$, $T(N/e)$, $T_C(N,e)$
and $T_L(N,e)$ we use, the resulting polynomial will belong to the
same equivalence class of ${\mathbb Z}[\Lambda]$ modulo $I_1$.
Thus we need only to show that one specific representative of
$T(M\otimes_+N)$ and $T(M\otimes_-N)$ can be computed from one
specific representative of $T(M)$, $T(N\backslash e)$, $T(N/e)$,
$T_C(N,e)$ and $T_L(N,e)$ by the above substitutions under one
special labelling rule. We will show this for $T(M\otimes_+N)$.
The proof for $T(M\otimes_-N)$ is exactly the same.

Assume that $M$ has $m$ edges and $N\B e$ has $n$ edges and both
have been assigned a labelling using positive integers from $1$ to
$m$ and $1$ to $n$ respectively. For the sake of convenience, an
edge of $M$ with label $i$ will be called $e_i$ in this proof. We
will now label $M\otimes_+N$ in the following way. First, any
negative edge in $M$ (which is not affected by the tensor
operation in $M\otimes_+N$) will retain its original label index.
On the other hand, after a positive edge with index $i$ is
replaced by $N\B e$, the edge in this copy of $N\B e$ with
original label index $j$ will now be assigned label index
$i+\frac{j}{n+1}$. For the sake of convenience, we may denote by
$N_i$ the copy of $N\B e$ that replaces the edge in $M$ with label
$i$ (as a part of the graph $M\otimes_+N$). Under this labelling,
if $e_j$ is a negative edge in $M$ (so it will not be replaced in
$T(M\otimes_+N)$ by $N\B e$), then the label of any edge in $N_i$
is less than $j$ if $i<j$, and is larger than $j$ if $i>j$. This
same rule applies to edges from $N_i$ and $N_j$ as well.

Let us now consider a spanning tree $P^\p$ of $M\otimes_+N$.
Notice that $P^\p$ induces a spanning tree $P$ of $M$ in a natural
way: if $P^\p$ contains an edge of $M$, then that edge is
retained. If $P^\p$ contains a path in $N_i$ that connects the two
vertices of the edge $e_i$ in $M$ replaced by $N_i$, then that
path is replaced by $e_i$, otherwise the edges of $N_i$ will
simply be removed. On the other hand, a spanning tree $P$ of $M$
can be extended to a spanning tree of $M\otimes_+N$ by replacing a
positive edge $e_i$ in $P$ by a spanning tree of $N$ (notice that
this spanning tree may contain the edge $e$ and in this case the
resulting spanning tree $P^\p$ of $M\otimes_+N$ would not contain
a path connecting the two vertices of the edge $e_i$). Any
spanning tree $P^\p$ of $M\otimes_+N$ obtained from the spanning
tree $P$ of $M$ is called an \textit{offspring} of $P$ and the
tree $P$ is called a \textit{parent tree} of $P^\p$. We will now
consider the total contribution of all the offsprings of a parent
tree $P$ in $M$.

Case 1. The negative edges of $M\otimes_+N$ that come directly
from $M$. Let $f$ be such an edge and let $P^\p$ be an offspring
of $P$. Apparently, $f$ is internal (external) to $P$ if and only
if it is internal (external) to $P^\p$. Furthermore, the activity
of $f$ in $P^\p$ is the same as its activity in $P$ by the choice
of our labelling for $M\otimes_+N$. Thus $f$ makes the same
contribution in $T(M)$ and $T(M\otimes_+N)$. In other words, we do
not replace the negatively indexed variables in $T(M)$ when
computing $T(M\otimes_+N)$.

Case 2. A positive edge $e_i$ in $M$ that is internally active
with respect to $P$. This edge is replaced by $N_i$. Notice that
any offspring $P^\p$ of $P$ induces a spanning tree $P_i$ of $N$
that does not contain $e$, that is, a spanning tree $P_i$ of $N\B
e$. In other words, $P^\p$ contains a path (in $N_i$) connecting
the two vertices of the edge $e_i$. Any internally active edge $f$
of $N_i$ with respect to $P_i$ is internal to $P^\p$ and is also
active since any edge $g$ such that $(P^\p \setminus  f ) \cup
 g $ is a spanning tree of $M\otimes_+N$ will have a larger
label by the labelling of $M\otimes_+N$. Of course, an internally
inactive edge $f$ of $N_i$ with respect to $P_i$ is an internally
inactive edge in $M\otimes_+N$ with respect to $P^\p$. On the
other hand, an externally active edge $f$ of $N_i$ with respect to
$P_i$ is external to $P^\p$ since $C(P^\p,f)=C(P_i,f)$ and the
labelling of $N_i$ inherits the relative ordering of the labelling
of $N\B e$ by our labelling choice. It follows that the activity
of $f$ in $N_i$ is the same as its activity in $M\otimes_+N$.
Therefore, the total contribution of $N_i$ to $T(M\otimes_+N)$
with respect to any given parent spanning tree $P$ is $T(N\B e)$.
That is, we may replace $x_+$ in $T(M)$ by $T(N\B e)$ when
computing $T(M\otimes_+N)$.

Case 3. A positive edge $e_i$ in $M$ that is internally inactive
with respect to $P$. Again such an edge is replaced by $N_i$ and
$P^\p$ contains a path (in $N_i$) connecting the two vertices of
the edge $e_i$. If $f$ is an edge on this path, then $f$ is
inactive since there exists an edge $g$ of $M\otimes_+N$ with
label smaller than that of $f$ such that $(P^\p \setminus  f )
\cup  g $ is a spanning tree of $M\otimes_+N$, see Figure
\ref{figconn}. The top of Figure \ref{figconn} shows the graph of
$P$, with $e_i \in P$, $e_j \notin P$, and $j < i$.  Thus $e_i$ is
internally inactive. The bottom left of Figure \ref{figconn} shows
the case when $e_j$ is negative and is not replaced by $N\B e$. In
this case, $(P^\p \setminus  f ) \cup  e_j $ is a spanning tree of
$M\otimes_+N$ with the label $e_j$ less than the label of $f$
(which is between $i$ and $i+1$ by our choice of labelling on
$M\otimes_+N$). The bottom right of Figure \ref{figconn} shows the
case when $e_j$ is positive and is replaced by $N_j=N\B e$. In
this case there exists some edge $f_j$ in $N_j$, $f_j \notin
P^\p$, such that $(P^\p \setminus  f ) \cup  f_j $ is a spanning
tree of $M\otimes_+N$, and again the label of $f_j$ (between $j$
and $j+1$) is less than the label of $f$ (between $i$ and $i+1$).
To summarize: an internally active edge $f$ in $N_j$ with respect
to $P_j$ is internally active in $M\otimes_+N$ with respect to
$P^\p$ if and only if it is not on $C(P_j,e)$. Finally, the
activity of any edge $f$ of $N_j$ external to $P^\p$ is decided
within $N_j$ since the unique cycle $C(P^\p,f)=C(P_j,f)$ is
contained in $N_j$. These rules are exactly the ones that define
the polynomial $T_L(N,e)$, hence we may replace $A_+$ in $T(M)$ by
$T_L(N,e)$ when computing $T(M\otimes_+N)$.

\begin{figure}[!htb]
\begin{center}
\includegraphics[scale=0.8]{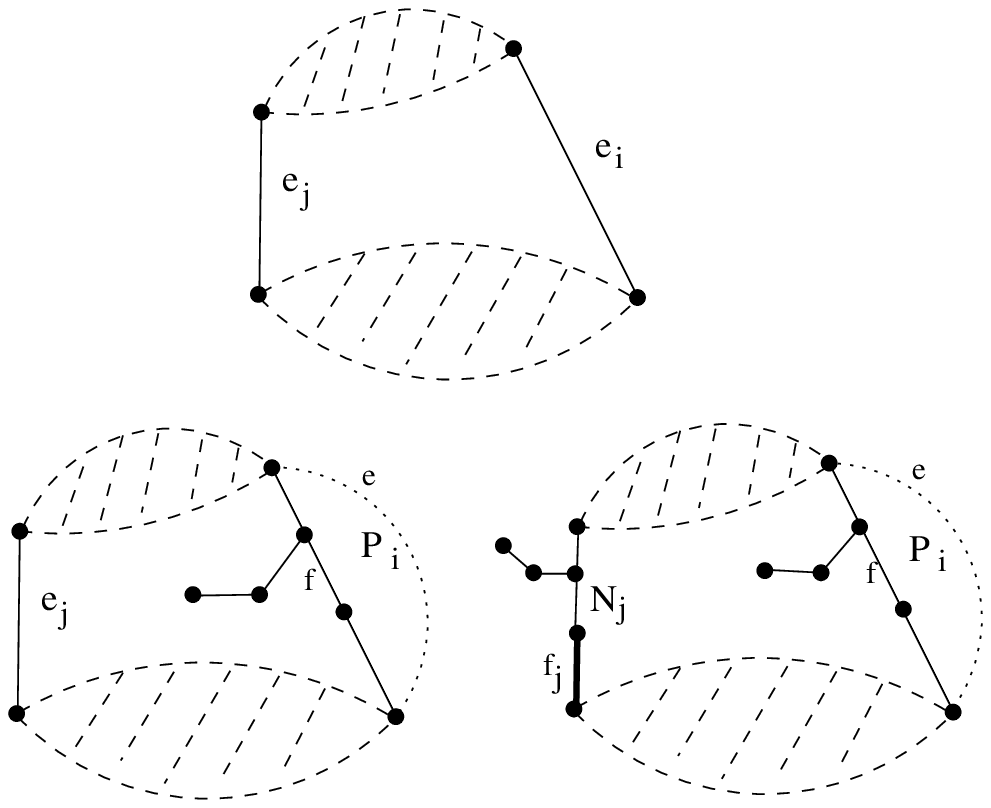}
\caption{The activities of an internal edge in $N_i$ when $e_i$ is
internally inactive in $M$.}\label{figconn}
\end{center}
\end{figure}

Case 4. A positive edge $e_i$ in $M$ that is externally active
with respect to $P$. Let $N_i$ be the copy of $N\B e$ replacing
$e_i$ with our labelling defined earlier and let $P^\p$ be an
offspring of $P$. Since $e_i$ is not on $P$, $P^\p$ does not
contain a path (in $N_i$) that connects the two vertices of $e_i$.
It follows that $P^\p$ induces a spanning tree $P_i^\p$ in $N/e$.
Let $f$ be an internal edge on $P_i^\p$. Then $P_i^\p \setminus
 f $ can be made into a tree only by adding an edge in $N_i$. In
other words, the activity of $f$ is decided ``locally" in $N/e$
with respect to the spanning tree $P^\p_i$. On the other hand, if
$f$ is an external edge in $M\otimes_+N$ with respect to $P^\p$,
then the cycle $C(P^\p,f)$ is contained in $N_i$ and can be
considered as a cycle of $N/e$. Thus the activity of $f$ in
$M\otimes_+N$ with respect to $P^\p$ is the same as the activity
of $f$ in $N/e$ with respect to $P^\p_i$. It follows that the
total contribution of the edges of $N_i$ (over all possible
$P^\p$) is the same as $T(N/e)$. Hence we may replace $y_+$ in
$T(M)$ by $T(N/e)$ when computing $T(M\otimes_+N)$.

Case 5. A positive edge $e_i$ in $M$ that is externally inactive
with respect to $P$. Let $N_i$ be the copy of $N\B e$ replacing
$e_i$ with our labelling defined earlier and let $P^\p$ be an
offspring of $P$. Again, $P^\p$ does not contain a path (in $N_i$)
that connects the two vertices of $e_i$ hence it induces a
spanning tree $P_i^\p$ in $N/e$. As in case 4 above, if $f$ is an
external edge in $M\otimes_+N$ with respect to $P^\p$, then the
cycle $C(P^\p,f)$ is contained in $N_i$ and can be considered as a
cycle of $N/e$. Thus the activity of $f$ in $M\otimes_+N$ with
respect to $P^\p$ is the same as the activity of $f$ in $N/e$ with
respect to $P^\p_i$. On the other hand, if $f$ is an internal edge
in $M\otimes_+N$ with respect to $P^\p$, then $f$ is internal in
$N/e$ with respect to $P_i^\p$. If $f$ is inactive in $N/e$
locally with respect to $P_i^\p$, then it is inactive. If the
addition of $f$ to $P^\p$ completes a path connecting the two
vertices of $e_i$, then there exists an edge in $M\otimes_+N$ (but
not in $N_i$) with a smaller label that is on the cycle
$C(P^\p,f)$. In any case, $f$ is still inactive in $M\otimes_+N$.
Finally, if $f$ is internally active in $N/e$ with respect to
$P_i^\p$, then it will still be internally active in $M\otimes_+N$
if the cycle $C(P^\p,f)$ is contained in $N_i$, otherwise $f$
would be inactive in $M\otimes_+N$. But this means exactly that
any externally active edge $f$ such that $e$ is on the cycle
$C(P^{\p\p}_i,f)$ in $N$ is to be considered as externally
inactive in $M\otimes_+N$ (here $P^{\p\p}_i$ is the spanning tree
of $N$ obtained from $P^\p_i$ by recovering the edge $e$). It
follows that the total contribution of the edges of $N_i$ (over
all possible $P^\p$) is the same as $T_C(N,e)$. Hence we may
replace $B_+$ in $T(M)$ by $T_C(N,e)$ when computing
$T(M\otimes_+N)$. This finishes our proof.
\end{proof}

We conclude this section by outlining how Theorem~\ref{Ttensor} and its
proof may be specialized to the unsigned case. Just like at the end of
Section~\ref{s_tcl}, consider unsigned graphs
as special signed graphs whose edges are all positive, and all
polynomial invariants as elements of the polynomial ring ${\mathbb
  Z}[x_+,y_+,A_+,B_+]$. Given an unsigned graph $M$ and an unsigned
graph $N$ with a distinguished edge $e$, considered as signed
graphs having positive edges only, the unsigned tensor product
$M\otimes N$ is identifiable with  the signed tensor product
$M\otimes_+N$. (For a definition of the unsigned tensor product
see~\cite{B}, \cite{BO}, \cite{Hu}, or \cite{Ja}, we will use the
above observation as our definition.) We may adapt the proof of
Theorem~\ref{Ttensor} in such a way that negative variables never
appear in the picture, and obtain the following analogous
statement.
\begin{theorem}\label{Tutensor}
Let $M$ be an unsigned graph and $N$ an unsigned graph with a
distinguished edge $e$. Then $T(M\otimes_+N)\in {\mathbb
Z}[x_+,y_+,A_+,B_+]$ can be computed from $T(M)\in {\mathbb
Z}[x_+,y_+,A_+,B_+]$ by using the substitutions
$$
x_+\mapsto T(N\backslash e)\quad
y_+\mapsto T(N/e)\quad
A_+\mapsto T_L(N,e)\quad
B_+\mapsto T_C(N,e).
$$
\end{theorem}
Observe now that the signed Tutte polynomial of $M$, when written
as a sum of monomials in positive variables, has the property that
the total degree in $x_+$ and $A_+$ of each monomial is the total
number of internal edges of some spanning tree, a.k.a. the {\em
rank} $r(M)$ of the matroid $M$, whereas the total degree in $y_+$
and $B_+$ is the total number of external edges, that is,
$|M|-r(M)$. After extending ${\mathbb Z}[x_+,y_+,A_+,B_+]$ to its
localization ${\mathbb Z}[x_+,y_+,A_+,B_+]_S$ by the semigroup $S$
generated by $\{T_C(N,e),T_L(N,e)\}$, Theorem~\ref{Tutensor} may
be rephrased as follows.
\begin{corollary}
\label{C_utensors}
Let $M$ be an unsigned graph and $N$ an unsigned graph with distinguished edge
$e$. Then $T(M\otimes_+N; x_+,y_+,A_+,B_+)\in {\mathbb Z}[x_+,y_+,A_+,B_+]_S$
is given by
$$
T(M\otimes_+N)=T_L(N,e)^{r(M)}T_C(N,e)^{|M|-r(M)}
\cdot T(M; T(N\backslash e)/T_L(N,e), T(N/e)/T_C(N,e), 1,1).
$$
Here all polynomials are considered as elements of ${\mathbb
  Z}[x_+,y_+,A_+,B_+]_S$.
\end{corollary}
It should be noted that for a nontrivial graph $N$ with a
distinguished edge $e$, the polynomials $T_C(N,e)$ and $T_C(N,e)$
are not zero since they may be computed as a sum of monomials of
positive coefficients. This new substitution rule has the property
that the contribution of the inactive edges of $M$ is $1$, hence
the formula ``factors'' through using the ordinary Tutte
polynomial of $M$. More precisely, we may uniquely extend the
homomorphism $\phi:{\mathbb
  Z}[x_+,y_+,A_+,B_+]\rightarrow {\mathbb Z}[x,y]$ given by $x_+\mapsto
x$, $y_+\mapsto y$, $A_+\mapsto 1$, $B_+\mapsto 1$ to a homomorphism ${\mathbb
  Z}[x_+,y_+,A_+,B_+]_S\rightarrow {\mathbb Z}[x,y]_{\phi(S)}$ and
apply it to both sides of the equation in
Corollary~\ref{C_utensors}. Observe that both sides belong to the
subring ${\mathbb Z}[x,y]$, only some of the calculation on the
right hand side needs to be performed in
${\mathbb Z}[x,y]_{\phi(S)}$. Thus we obtain a new
proof of  the following classical result:

\begin{corollary}
\label{Ctensor} Let $M$ be an unsigned graph and $N$ an unsigned
graph with a distinguished edge $e$. Then the ordinary Tutte
polynomial $T(M\otimes N)\in {\mathbb Z}[x,y]$ may be obtained
from the ordinary Tutte polynomial $T(M)\in {\mathbb Z}[x,y]$ by
substituting $T(N\backslash e)/T_L(N,e)$ into $x$,
$T(N/e)/T_C(N,e)$ into $y$, and multiplying the resulting rational
expression with $T_L(N,e)^{r(M)}T_C(N,e)^{|M|-r(M)}$. Here
$T_C(N,e)$ and $T_L(N,e)$ are elements of ${\mathbb Z}[x,y]$,
defined by (\ref{E_utctl}).
\end{corollary}
This result may be found in~\cite[Equation (4.1)]{Ja} without a
proof. The first published proof of a differently phrased but
equivalent statement is due to Brylawski~\cite{B}, which appears
to take a very different approach. It should be noted that we
referred to our signed results only for brevity's sake, the
shortest direct proof of Corollary~\ref{Ctensor} may be obtained
by directly adapting the proof of Theorem~\ref{Ttensor} to the
unsigned case. We leave the details to the reader.

\section{Applications to knot theory}
\label{s_app}

\begin{figure}[!htb]
\begin{center}
\includegraphics[scale=0.6]{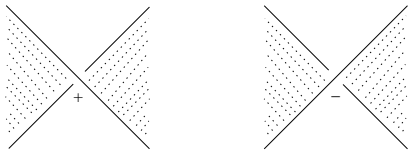}
\caption{The assignment of signs at a crossing (vertex) for the
graph $G$.}\label{crossingsign}
\end{center}
\end{figure}

It is well known that the Jones polynomial of a knot $K$ can be
obtained from the Tutte polynomial of the dual graph of $D$, where
$D$ is a regular projection of $K$ \cite{K2}. Let us first give a
brief description of this process. One starts from a regular
projection $D$ of the knot $K$. We then shade the regions in its
projection either ``white'' or ``dark'' in a checkerboard fashion,
so that no two dark regions are adjacent, and no two white regions
are adjacent. We usually consider the infinite region surrounding
the knot projection to be white. Note that as we move diagonally
over a knot crossing, we go from a white region to a white region,
or from a dark region to a dark region.  Next we construct a dual
graph of $D$ by converting the dark regions in $D$ into vertices
in a graph $G$ and converting the crossings in $D$ between two dark
regions into edges incident to the corresponding vertices in $G$.
So if we can move diagonally over a knot crossing from one dark
region to another, then these two dark regions and the crossing
will be represented in $G$ as two vertices connected by an edge.
Note that we may obtain parallel edges from some knot projections.
Now we have our unsigned graph. To obtain the signed version, we
look at each crossing in the knot projection.  If, after the upper
strand passes over the lower, the dark region is to the left of
the upper strand, then we denote this as a positive crossing.  If
the dark region is to the right of the upper strand, we denote it
as a negative crossing.  Then our signed graph is obtained by
marking each edge of $G$ with the same sign as the crossing of $K$
to which it corresponds. See Figure \ref{crossingsign}.

The following theorem is due to Kauffman \cite{K1,K2}.

\begin{theorem}\label{T41}
Let $G$ be the (signed) dual graph of a regular knot projection
$D$ of $K$ as described above, then $T(G)$ equals the Kauffman
bracket polynomial $\langle K\rangle$ under the following variable
substitutions:
\begin{eqnarray*}
&&x_+\mapsto -A^{-3},\ x_-\mapsto -A^3,\
y_+\mapsto -A^3,\ y_-\mapsto-A^{-3}\\
&&A_+\mapsto A,\ A_-\mapsto A^{-1},\ B_+\mapsto A^{-1},\
B_-\mapsto A.
\end{eqnarray*}
Furthermore, the Jones polynomial $V_K(t)$ of $K$ can be obtained
from
\begin{equation}\label{eq51}
V_K(t)=(-A^{-3})^{w(K)}\langle K\rangle
\end{equation}
by setting $A=t^{-\frac{1}{4}}$, where $w(K)$ is the writhe of the
projection $D$.
\end{theorem}

It is thus possible for us to use Theorem \ref{Ttensor} to compute
the Jones polynomials for some large non-alternating knots. We
will demonstrate this by a few examples. We will first do this for
a small knot so we can compare our result with the direct
computation result using an existing software. We will then do
this for a much larger knot beyond the capacity of the existing
programs.

\begin{example}\em
In this example, let us use the signed graphs obtained from the
standard minimal knot diagrams of the non-alternating knot
$9_{49}$ and the alternating knot $4_1$ as $M$ and $N$
respectively, see Figure \ref{945}. Notice that the identifying
edge $e$ is marked in the figure.
\begin{figure}[!htb]
\begin{center}
\includegraphics[scale=0.8]{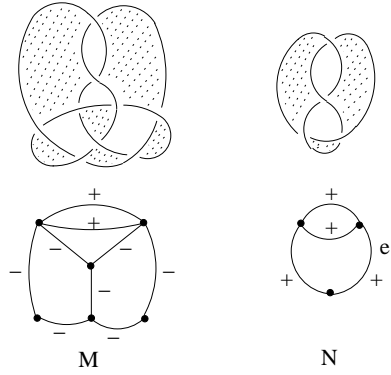}
\caption{The knots $9_{49}$, $4_1$ and their corresponding signed
graphs.}\label{945}
\end{center}
\end{figure}
Since the tensor product is not unique, it is possible for us to
get several different knots from $M\otimes_+N$. One of such is
shown in Figure \ref{mn41}.
\begin{figure}[!htb]
\begin{center}
\includegraphics[scale=1.0]{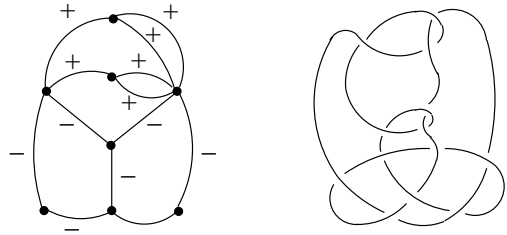}
\caption{A knot obtained from positive tensor product using
$9_{49}$ and $4_1$.}\label{mn41}
\end{center}
\end{figure}
The Tutte polynomial for the signed graph corresponding to $M$ in
this case is (using its $55$ spanning trees):
\begin{eqnarray*}
& &B_+^2(2x_-^2A_-^3y_-B_- + 3x_-^3A_-^2B_-^2 + x_-^5B_-^2 +
2x_-^4A_-B_-^2\\
& +& 2x_-A_-^4y_-B_- + 2x_-^2A_-^3B_-^2 + x_-A_-^4B_-^2 +
A_-^5y_-^2 + A_-^5y_-B_-)\\
& +& (A_+y_+ + A_+B_+)(2x_-^3A_-B_-^3 + 3x_-^2A_-^2B_-^3 +
4x_-A_-^3y_-B_-^2\\
& +& x_-A_-^3B_-^3 + 2A_-^4y_-^2B_- + A_-^4y_-B_-^2 + x_-^4B_-^3\\
&+& x_-^3A_-y_-B_-^2 + 2x_-^2A_-^2y_-B_-^2 + A_-^4y_-^3 +
2x_-A_-^3y_-^2B_-).
\end{eqnarray*}
To verify that this is correct, we calculate the bracket
polynomial from it using Theorem \ref{T41}. This gives us
$$
-A^{19} + 2A^{15} - 4A^{11} + 4A^7 - 5A^3 + 4A^{-1} - 3A^{-5} +
2A^{-9}.
$$
Since the writhe of this projection of $9_{49}$ is 9, we obtain
the Jones polynomial of $9_{49}$:
\begin{eqnarray*}
V(9_{49})&=&(-A^{-3})^{w(9_{49})}\langle 9_{49} \rangle \\
&=& -A^{-27}(-A^{19} + 2A^{15} - 4A^{11} + 4A^7 - 5A^3 + 4A^{-1} -
3A^{-5} + 2A^{-9})\\
&=& A^{-8} - 2A^{-12} + 4A^{-16} - 4A^{-20} + 5A^{-24} - 4A^{-28}
+ 3A^{-32} - 2A^{-36}\\
& =& t^2 - 2t^3 + 4t^4 - 4t^5 + 5t^6 - 4t^7 + 3t^8 - 2t^9.
\end{eqnarray*}
This matches the known Jones polynomial for $9_{49}$.

We next find the polynomials $T(N\setminus e)$, $T(N/e)$,
$T_L(N,e)$ and $T_C(N,e)$.  We calculate the first two by drawing
the spanning trees of $N\setminus e$ and $N/e$, then the latter
two are obtained from the first two by applying the exception
rules.  We have
\begin{eqnarray*}
T(N\setminus e) &=& x_+^2B_+ + x_+y_+A_+,\\
T(N/e) &=& x_+B_+^2 + y_+A_+B_+ + y_+^2A_+, \\
T_L(N,e) &=& A_+^2B_+ + A_+^2y_+,\\
T_C(N,e) &=& x_+B_+^2 + A_+B_+^2 + y_+A_+B_+.
\end{eqnarray*}
We may now replace $B_+$, $y_+$ and $A_+$ in $T(M)$ with
$T_C(N,e)$, $T(N/e)$, and $T_L(N,e)$ (since $x_+$ does not appear
in $T(M)$) respectively to obtain $T(M \otimes_+ N)$:
\begin{eqnarray*}
T(M \otimes_+ N)& =& (x_+B_+^2 + A_+B_+^2 +
y_+A_+B_+)^2(2x_-^2A_-^3y_-B_- \\
&+& 3x_-^3A_-^2B_-^2 + x_-^5B_-^2 + 2x_-^4A_-B_-^2 +
2x_-A_-^4y_-B_- \\
&+& 2x_-^2A_-^3B_-^2 + x_-A_-^4B_-^2 + A_-^5y_-^2 + A_-^5y_-B_-)
\\
&+& [(A_+^2B_+ + A_+^2y_+)(x_+B_+^2 + y_+A_+B_+ + y_+^2A_+) \\
&+& (A_+^2B_+ + A_+^2y_+)(x_+B_+^2 + A_+B_+^2 +
y_+A_+B_+)](2x_-^3A_-B_-^3 \\
&+& 3x_-^2A_-^2B_-^3 + 4x_-A_-^3y_-B_-^2 + x_-A_-^3B_-^3 +
2A_-^4y_-^2B_- \\
&+& A_-^4y_-B_-^2 + x_-^4B_-^3 + x_-^3A_-y_-B_-^2 +
2x_-^2A_-^2y_-B_-^2 \\
&+& A_-^4y_-^3 + 2x_-A_-^3y_-^2B_-).
\end{eqnarray*}
We can then calculate the bracket polynomial and the Jones
polynomial for $K_{M \otimes_+ N}$:
\begin{eqnarray*}
\langle K_{M \otimes_+ N} \rangle & = & -A^{27} + 4A^{23} -
9A^{19} + 14A^{15} - 17A^{11} + 19A^7\\
& -& 18A^3 + 13A^{-1} - 9A^{-5} + 3A^{-9} - A^{-17} + A^{-21}.
\end{eqnarray*}
The writhe of the knot projection diagram corresponding to ${M
\otimes_+ N}$ is 1. Hence
\begin{eqnarray*}
&& V_{K_{M \otimes_+ N}}(t) =(-A^{-3})^{w(K_{M \otimes_+ N})}
\langle K_{M \otimes_+ N} \rangle\\
& =& A^{24} - 4A^{20} + 9A^{16} - 14A^{12} + 17A^8 - 19A^4 + 18 -
13A^{-4} + 9A^{-8} - 3A^{-12} + A^{-20} - A^{-24} \\
&=& t^{-6} - 4t^{-5} + 9t^{-4} - 14t^{-3} + 17t^{-2} - 19t^{-1} +
18 - 13t + 9t^2 - 3t^3 + t^5 - t^6.
\end{eqnarray*}
Using the DT code $[-4,-22,-8,-20,14,24,18,26,10,-6,-2,12,16]$
obtained from the diagram in Figure \ref{mn41}, we get the same
Jones polynomial using Knotscape \cite{knotscape}.
\end{example}

\begin{example}\em
Our next example deals with a knot family with a parameter $k$. We
will illustrate the details using the case of $k=3$ (illustrated
in Figure \ref{fig5}) and then show our computation results for
the cases of $k=5$, $7$ and $9$. Notice that the projection
diagram $D$ is non-alternating. For $k=3$, there are 19 crossings
in the diagram, 10 positive and 9 negative with respect to the
shaded regions. The corresponding signed graph $G$ of this diagram
(using the shaded regions in the figure as the vertices) is shown
on the right side. A labelling of the edges of $G$ is given in the
figure as well.

\begin{figure}[!htb]
\begin{center}
\includegraphics[scale=0.8]{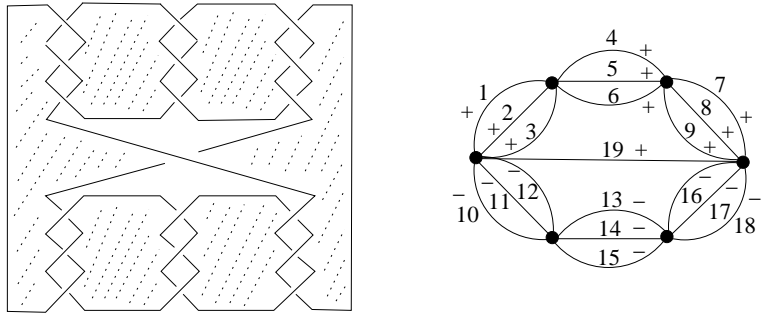}
\caption{A 19 crossing knot diagram and its corresponding signed
graph $G$.}\label{fig5}
\end{center}
\end{figure}

\begin{figure}[!htb]
\begin{center}
\includegraphics[scale=0.8]{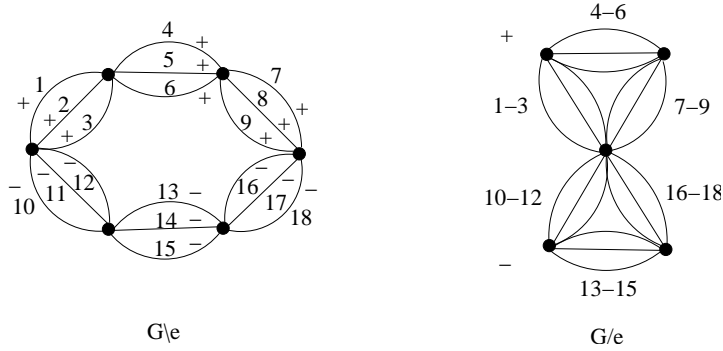}
\caption{The graphs obtained from the graph $G$ in Figure
\ref{fig5} by deletion and contraction.}\label{fig6}
\end{center}
\end{figure}

Denote the edge labelled by 19 in the figure by $g$. Let
$G\backslash g$ be the graph obtained from $G$ by deleting $g$ and
let $G/g$ be the graph obtained from $G$ by contracting $g$, that
is, the vertices incident to $g$ are identified as one single
vertex after $g$ is removed. See Figure \ref{fig6}. By the
recursive formula (\ref{recur}), we have
$$
T(G)=B_+ T(G\backslash g)+A_+T(G/g).
$$
Thus we will concentrate on computing $T(G\backslash g)$ and
$T(G/g)$. Notice that $G\backslash g$ and $G/g$ can be obtained
from the simple graphs $M_1$ and $M_2$ shown in Figure \ref{m12}
by applying repeated tensor product operations: $ G\backslash
g=((((M_1\otimes_+S_+)\otimes_-S_-)\otimes_+T_+)\otimes_-T_- $ and
$ G/g=((((M_2\otimes_+S_+)\otimes_-S_-)\otimes_+T_+)\otimes_-T_-,
$ where $S_-$, $S_+$, $T_-$ and $T_+$ are shown in Figure \ref{ST}
for the case of $k=3$. Notice that the identifying edge $e$ is
marked in the figure.

\begin{figure}[!htb]
\begin{center}
\includegraphics[scale=0.8]{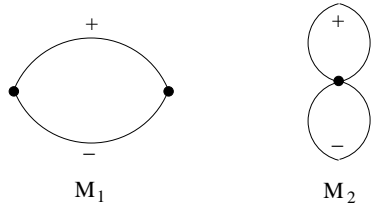}
\caption{Two very simple graphs.}\label{m12}
\end{center}
\end{figure}
\begin{figure}[!htb]
\begin{center}
\includegraphics[scale=0.8]{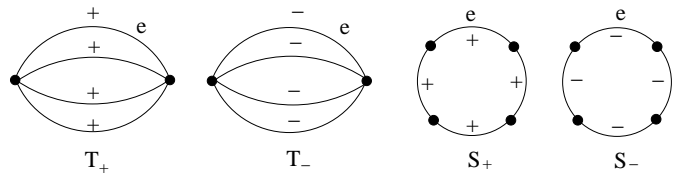}
\caption{These graphs are also called {\em stretching and
thickening} graphs.}\label{ST}
\end{center}
\end{figure}

Since $T(M_1)=x_+B_-+y_+A_-$, $T(M_2)=y_+y_-$ and each of $S_-$,
$S_+$, $T_-$ and $T_+$ uses only one sign, $T(G\backslash g)$ and
$T(G/g)$ can be computed by making the following substitutions
twice starting from $x_+B_-+y_+A_-$ and $y_+y_-$. (We leave it as
an exercise to our reader to compute the polynomials
$T(S_\pm\backslash e)$, $T(S_\pm/e)$, $T(T_\pm\backslash e)$,
$T(T_\pm/e)$, $T_C(S_\pm,e)$, $T_C(T_\pm,e)$, $T_L(S_\pm,e)$ and
$T_L(T_\pm,e)$.)
\begin{eqnarray*}
A_-&\mapsto & A_-(y_{-}^2+y_-B_- + B_{-}^2),\quad B_-\mapsto  B_-^3,\\
A_+&\mapsto & A_+(y_{+}^2+y_+B_+ + B_{+}^2),\quad B_+\mapsto  B_+^3,\\
x_-&\mapsto & B_{-}^{2}x_-+ (y_{-} + B_{-}) A_{-}y_{-},\quad y_-\mapsto y_-^3,\\
x_+&\mapsto & B_{+}^{2}x_{+} + (y_{+} + B_{+}) A_{+}y_{+},\quad
y_+\mapsto  y_+^3.
\end{eqnarray*}
After making the substitutions listed in Theorem \ref{T41} and
substituting $t^{-\frac{1}{4}}$ for $A$ in Equation (\ref{eq51})
(with $w(K)=-1$), we obtain
\begin{eqnarray*}
V_K(t)&=&t^{-10}(1-4t+12t^2-26t^3+49t^4-74t^5+96t^6-112t^7+110t^8\\
&-&97t^9+77t^{10}-47t^{11}+23t^{12}-8t^{13}-2t^{14}+3t^{15}-t^{16}+
t^{17}).
\end{eqnarray*}
This matches the computation result obtained by using
\textit{Knotscape} \cite{knotscape}. Use the same approach, we can
compute the Jones polynomials for larger values of $k$. The result
for $k=5$, $7$ and $9$ are listed below. The computation time for
the case of $k=9$ is about 10 minutes using Maple on a PC.

For $k=5$,
\begin{eqnarray*}
V_K(t)&=&t^{-26}(1 - 6t + 26t^2 - 91t^3 + 275t^4 - 737t^5 +
1796t^6 -
4021t^7 + 8366t^8 - 16284t^9\\
&+& 29818t^{10} - 51606t^{11} + 84676t^{12} - 132106t^{13} +
196368t^{14} - 278544t^{15} + 377546t^{16}\\
&-& 489336t^{17} + 606846t^{18} - 720177t^{19} + 817720t^{20} -
887911t^{21} + 920952t^{22} - 911068t^{23}\\
&+& 857489t^{24} - 765053t^{25} + 643579t^{26} - 505933t^{27} +
366267t^{28} - 237242t^{29} + 128459t^{30}\\
&-& 45354t^{31} - 11121t^{32} + 43431t^{33} - 56574t^{34} +
56418t^{35}
- 48576t^{36} + 37646t^{37}\\
&-& 26696t^{38} + 17478t^{39} - 10594t^{40} + 5941t^{41} -
3081t^{42} +
1466t^{43} - 637t^{44}\\
&+& 250t^{45} - 86t^{46} + 26t^{47} - 6t^{48} + t^{49}).
\end{eqnarray*}

For $k=7$,
\begin{eqnarray*}
V_K(t)&=&t^{-50}(1 - 8t + 43t^2 - 183t^3 + 666t^4 - 2157t^5 +
6370t^6 -
17425t^7 + 44654t^8 - 108067t^9\\
&+& 248536t^{10} - 545847t^{11} + 1149387t^{12} - 2328122t^{13} +
4548764t^{14} - 8593271t^{15}\\
&+& 15728483t^{16} - 27941544t^{17} + 48253003t^{18} -
81115378t^{19} + 132896097t^{20} - 212430488t^{21}\\
\phantom{V_K(t)}&+&331612373t^{22} - 505966329t^{23} +
755122019t^{24} -
1103084529t^{25} + 1578177868t^{26}\\
&-& 2212528476t^{27} + 3040964638t^{28} - 4099238067t^{29} +
5421525110t^{30} - 7037249318t^{31}\\
&+& 8967357925t^{32} - 11220302612t^{33} + 13788073932t^{34} -
16642729060t^{35} + 19733901580t^{36}\\
&-& 22987768175t^{37} + 26307888783t^{38} - 29578176531t^{39} +
32668072879t^{40} - 35439739995t^{41} \\
&+& 37756834115t^{42} - 39494183464t^{43} + 40547500331t^{44} -
40842187270t^{45} + 40340281916t^{46}\\
&-& 39044758086t^{47} + 37000613964t^{48} - 34292527561t^{49} +
31039223375t^{50} - 27385045791t^{51}\\
&+& 23489568932t^{52} - 19516242140t^{53} + 15621197678t^{54} -
11943255900t^{55} + 8596013866t^{56}\\
&-& 5662618702t^{57} + 3193491734t^{58} - 1206952613t^{59} -
307642091t^{60} + 1384824407t^{61}\\
&-& 2076750500t^{62} + 2446499548t^{63} - 2561637408t^{64} +
2488610216t^{65} - 2288339306t^{66}\\
 &+& 2013210200t^{67} -
1705443587t^{68} + 1396707072t^{69} -
1108711902t^{70} + 854497057t^{71}
\\
 &-& 640099923t^{72} +
466343256t^{73} - 330533561t^{74} + 227925671t^{75} -
152882699t^{76}\\
 &+& 99711533t^{77} - 63198913t^{78} +
38898252t^{79} - 23227709t^{80} + 13441772t^{81} -
7528452t^{82}\end{eqnarray*}
\begin{eqnarray*}
&+& 4074514t^{83} - 2126998t^{84} + 1068635t^{85} - 515390t^{86} +
237854t^{87} - 104637t^{88} + 43667t^{89}\\
&-& 17180t^{90} + 6321t^{91} - 2150t^{92} + 666t^{93} - 183t^{94}
+ 43t^{95} - 8t^{96} + t^{97}).
\end{eqnarray*}

For $k=9$, the polynomial is too large to list, so we only list a
few terms below:
\begin{eqnarray*}
V_K(t)&=&t^{-82}(1 - 10t + 64t^2 - 319t^3 + 1345t^4 - 5008t^5 +
    \cdots \\
& & +\cdots  - 20193935024459t^{97} -
    101497138129454t^{98} + \cdots\\
& & \cdots - 5008t^{156} +
    1345t^{157} - 319t^{158} + 64t^{159} - 10t^{160} + t^{161}).
\end{eqnarray*}
\end{example}

The knots constructed this way are non-alternating since the Jones
polynomials are not alternating. It follows that the crossing
numbers of these knots (at least for the ones we have computed
above) are at least $2k^2$ since the breadth of the polynomials is
$2k^2-1$ and the crossing number of a non-alternating knot is
strictly larger than the breadth of its Jones polynomial.

The computation of the above example with $k=9$ took only minutes
on a PC using Maple. However, it is not always possible to compute
the Tutte polynomial quickly through the variable substitutions
given in the last section. It would be interesting to see what
(non-alternating) knots can be constructed using graph tensor
product whose Jones polynomials (or the breadths of their Jones
polynomials) can be computed with polynomial runtime. This is a
possible direction of future research. The authors also intend to
investigate the possibilities of extending their results to other
knot polynomials or other graph invariants that are related to the
Tutte polynomial.

\section*{Acknowledgement}
We wish to thank Professor T.\ G.\ Lucas for providing the main
idea behind Lemma~\ref{L_ip}.

\end{document}